\documentclass[openany, amssymb, psamsfonts]{amsart}

\usepackage[utf8]{inputenc}
\usepackage[margin=3cm,top=2.5cm,headheight=16pt,headsep=0.1in,heightrounded]{geometry}
\usepackage{hyperref}
\usepackage{physics}
\usepackage[most]{tcolorbox}
\usepackage{xparse}
\usepackage{accents}
\usepackage{pdfpages}
\usepackage{wrapfig}
\usepackage{changepage}
\usepackage{amsmath, amssymb, amsthm, amsfonts}
\usepackage{mathrsfs}
\usepackage{xcolor}
\usepackage{bbm}
\usepackage{tikz-cd}
\usepackage{tikz}
\usepackage{quiver}
\usepackage{mathtools}

\newtheorem{theorem}{Theorem}[section]
\newtheorem{definition}{Definition}[section]
\newtheorem{lemma}{Lemma}

\newcommand{\bb}[1]{\mathbb{#1}}
\newcommand{\N}{\bb{N}}
\newcommand{\Z}{\bb{Z}}

\newcommand{\s}{\bb{S}}
\newcommand{\R}{\mathbb{R}}
\newcommand{\C}{\bb{C}}

\newcommand{\D}{\Delta}
\newcommand{\ep}{\varepsilon}
\newcommand{\ts}[1]{\textbf{#1}}
\newcommand{\p}{\partial}

\hypersetup{%
  bookmarksnumbered=true,%
  bookmarks=true,%
  colorlinks=true,%
  linkcolor=black,%
  citecolor=purple,%
  filecolor=blue,%
  menucolor=blue,%
  pagecolor=blue,%
  urlcolor=blue,%
  pdfnewwindow=true,%
  pdfstartview=FitBH}

\title{Revisiting Pontryagin's Proof of Stable Stems 1 and 2}
\author{{\textsc{Trishan Mondal}}}
\date{}

\begin{document}
\maketitle
\newcommand{\ocr}{\mathscr{O}}

\begin{abstract}
    In this paper, we introduce fundamental notions of homotopy theory, including homotopy excision and the Freudenthal suspension theorem. We then explore framed cobordism and its connection to stable homotopy groups of spheres through the Pontryagin-Thom construction. Using this framework, we compute the stable stems in dimensions $0$, $1$, and $2$. This work is primarily expository, revisiting proofs from \cite{Pont1} with slight modifications incorporating modern notation. Furthermore, in the final section, we discuss 2-dimensional framed manifolds with Arf invariant one and examine why the result of \cite{Pont2} regarding $\pi_2^S$ is incorrect.
\end{abstract}

\tableofcontents

\section{\hspace{0.1cm} Introduction to Homotopy theory}
\noindent We begin by introducing the notion of higher relative homotopy groups. In this paper the category of topological spaces is denoted by \(\mathbf{Top}\), while \(\mathbf{hTop}\) represents its homotopy category, where morphisms are continuous maps considered up to homotopy equivalence. Similarly, \(\mathbf{Top}_{\ast}\) and \(\mathbf{hTop}_{\ast}\) denote the categories of based spaces and their corresponding homotopy categories, respectively. \begin{itemize}
    \item[] The \textbf{Homotopy group} $\pi_n : \mathbf{Top}_{\ast} \to \textbf{Groups}$ is a functor defined by $\pi_n(X,x_0) = [(\s^n,e),(X,x_0)]$. Where $[(\s^n,e),(X,x_0)]$ means the collection of maps from $\s^n \to X$ so that $f(e)=x_0$ upto homotopy equivalence. Here the group operation $[f]+[g]$ is given by homotopy class of $[(f\vee g )\circ c]$, \[
     \includegraphics[width = 7cm]{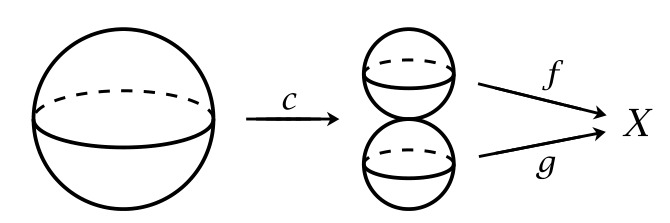}    
    \] here $c$ is the pinching map, pinched the equator to get $(\s^n)\vee (\s^n)$.

    \item If the space $X$ is path connected the homotopy group is independent of the base point $x_0$.
    \item For $n\geq 2$, homotopy groups are abelian. This follows from Eckmann-Hilton argument [proved \href{trishan8.github.io/resources/semsters/algtop2.pdf}{\textcolor{purple}{here}}]. 
    \item If $p: (\tilde{X},\tilde{x}_0)\to (X,x_0)$ is a covering then it induces isomorphism on homotopy groups for $n \geq 2$.
    \item Homotopy groups commutes with product, i.e. $\pi_n(X\times Y,(x_0,y_0)) \simeq \pi_n(X,x_0)\times \pi_n(Y,y_0)$.
\end{itemize}

\subsection{Relative Homotopy groups.}

 \noindent Given two spaces $X,Y \in \textbf{Top}_{\ast}$ and a  map $f : X\to Y$ (must be a based map), we can define \textbf{homotopy fiber} to be the pullback of the following diagram,  \[\begin{tikzcd}
	X & Y \\
	Ff & PY
	\arrow["f", from=1-1, to=1-2]
	\arrow["p"', from=2-2, to=1-2]
	\arrow[dashed, from=2-1, to=2-2]
	\arrow["\pi",dashed, from=2-1, to=1-1]
\end{tikzcd}\] Here $PY$ is the path space. Since, $Y$ is a based space at a point say $y_0$, $PY:= \qty{\gamma:([0,1],0) \to (Y,y_0)}$. The map $p : PY \to Y$ is given by $\gamma \mapsto \gamma(1)$. In other words $Ff$ is fiber product of $X$ and $PY$. We know the continuous map $f$ can be homotped to a fibration $g: E\to Y$ then the fibre of this fibration is homotopic to $Ff$. In the pullback diagram the map $\pi$ is given by the projection $(x,\gamma)\mapsto x$. 

\vspace*{0.2cm}

\noindent Now we define $\Omega (X,x_0) := \qty{\gamma \in PX : \gamma(1)=x_0}$ the loop space of the based space $X$. For any map $f: X \to Y$ there is a sequence of space as follows $$\cdots \Omega^2X \xrightarrow{\Omega^2 f}\Omega^2Y \xrightarrow{-\Omega i} \Omega Ff \xrightarrow{-\Omega\pi} \Omega X \xrightarrow{-\Omega f}\Omega Y\xrightarrow{i}Ff \xrightarrow{\pi}X \xrightarrow{f}Y $$ here $i$ is the natural inclusion $\gamma \mapsto (y_0,\gamma)$. In the above sequence any three consecutive spaces are part of fibration (upto homotopy). \textbf{Note:}  In general given any fibration $F \hookrightarrow E \xrightarrow{p} B$.

\begin{theorem}\label{thm:fibseq}
    For any space $S \in \mathbf{Top}_{\ast}$ the above fiber sequence induces the following long exact sequence,
    $$\cdots \to [S,\Omega Ff] \to [S,\Omega X]\to [S,\Omega Y] \to [S,Ff] \to [S,X]\to [S,Y]$$ 
\end{theorem}

\noindent Now we will define a functor $\Sigma : \textbf{Top}_{\ast} \to \textbf{Top}_{\ast}$ (called reduced suspension functor). If $X,Y$ are two based space we can define the smash product $X \wedge Y = X \times Y / (X \vee Y)$. For any space $X$, we define $\Sigma X = X \wedge S^1$. It's not hard to see as a functor $\Sigma$ and $\Omega$ are adjoint, i.e $$[\Sigma X, A]= [X, \Omega A]$$ for any map $f : X \to \Sigma A$, $f(x)$ is a loop based at $a_0$, we can define $(x,t)\mapsto f(x)(t)$ and this gives us a map from $\Sigma X \to A$, similarly any map $f: \Sigma X \to A$ will give us a map $x \mapsto f(x,t)$ (here $t$ varies over $S^1$ to give us the loop). This is the idea to establish the adjoint property. Recall, $\Sigma \s^n \simeq \s^{n+1}$. So from the definition of homotopy groups and the loop-suspension adjunction it follows $\pi_n(X,x_0)\simeq \pi_{n-1}(\Omega(X,x_0))$.

\vspace*{0.2cm}

\vspace*{0.2cm}

\noindent Consider $i$ to be inclusion of $x_0 \to X$. From the desciption of fibre product/homotopy fiber $Fi$ we get, $Fi = \qty{\gamma \in PX : \gamma(1)=x_0} = \Omega X$. We have also seen $\pi_n(X,x_0)= \pi_{n-1}(\Omega(X,x_0))$. This motivates us to give the definition of relative homotopy groups $\pi_n(X,A)$. Let, $i : (A,x_0) \hookrightarrow X$ then define \textbf{relative homotopy groups} (for $n \geq 1$)$$\pi_n(X,A):=\pi_{n-1}(Fi)$$From \ref{thm:fibseq} we can say there is the following long exact sequence, 

\[\begin{tikzcd}
	\cdots & {[\s^1,\Omega Fi]} & {[\s^1,\Omega A]} & {[\s^1,\Omega X]} & {[\s^1,Fi]} & {[\s^1,A]} & {[\s^1,X]} \\
	\cdots & {\pi_1(\Omega Fi)} & {\pi_1(\Omega A)} & {\pi_1(\Omega X)} & {\pi_1(Fi)} & {\pi_1(A)} & {\pi_1(X)} \\
	\cdots & {\pi_3(X,A)} & {\pi_2(A)} & {\pi_2(X)} & {\pi_2(X,A)} & {\pi_1(A)} & {\pi_1(X)}
	\arrow[from=1-2, to=1-3]
	\arrow["\simeq"', color={rgb,255:red,153;green,92;blue,214}, from=1-2, to=2-2]
	\arrow[from=1-3, to=1-4]
	\arrow["\simeq", color={rgb,255:red,153;green,92;blue,214}, from=1-3, to=2-3]
	\arrow[from=1-4, to=1-5]
	\arrow["\simeq", color={rgb,255:red,153;green,92;blue,214}, from=1-4, to=2-4]
	\arrow[from=1-5, to=1-6]
	\arrow["{\simeq }"', color={rgb,255:red,153;green,92;blue,214}, from=1-5, to=2-5]
	\arrow[from=1-6, to=1-7]
	\arrow["\simeq", color={rgb,255:red,153;green,92;blue,214}, from=1-6, to=2-6]
	\arrow["\simeq"', color={rgb,255:red,153;green,92;blue,214}, from=1-7, to=2-7]
	\arrow[from=2-2, to=2-3]
	\arrow["\simeq"', color={rgb,255:red,214;green,92;blue,92}, from=2-2, to=3-2]
	\arrow[from=2-3, to=2-4]
	\arrow["\simeq"', color={rgb,255:red,214;green,92;blue,92}, from=2-3, to=3-3]
	\arrow[from=2-4, to=2-5]
	\arrow["\simeq"', color={rgb,255:red,214;green,92;blue,92}, from=2-4, to=3-4]
	\arrow[from=2-5, to=2-6]
	\arrow["\simeq"', color={rgb,255:red,214;green,92;blue,92}, from=2-5, to=3-5]
	\arrow[from=2-6, to=2-7]
	\arrow["\simeq", color={rgb,255:red,214;green,92;blue,92}, from=2-6, to=3-6]
	\arrow["\simeq"', color={rgb,255:red,214;green,92;blue,92}, from=2-7, to=3-7]
	\arrow[from=3-2, to=3-3]
	\arrow[from=3-3, to=3-4]
	\arrow[from=3-4, to=3-5]
	\arrow[from=3-5, to=3-6]
	\arrow[from=3-6, to=3-7]
\end{tikzcd}\]

\noindent We can summarize the above discussion with the following theorem,

\begin{theorem} \label{thm:lespair}
 For a pair $(X,A) \in \textbf{Top}_{\ast}^2$, we have the following long exact sequence of homotopy groups
    \[
\cdots \to \pi_n(A) \to \pi_n(X) \to \pi_n(X,A)\xrightarrow{\p_{\ast}} \pi_{n-1}(A)\to \cdots    
\]
Infact for a fibration $E\to B$ with fiber $F$ then we have a Long exact sequence of homotopy groups, $$\cdots \to \pi_n(F)\to \pi_n(E)\to \pi_n(B)\xrightarrow{\p_{\ast}} \pi_{n-1}(F)\to \cdots$$
\end{theorem}

\noindent \textbf{Description of $\p_{\ast}$:} It's not hard to see this definition above is equivalent to defining $\pi_n(X,A) = [(I^n,\p I^n,J^n),(X,A,x_0)]$ where $J^n = \overline{I^n\setminus I^{n-1}\times\qty{1}}$. Any map $f: (I^n,\p I^n,J^n) \to (X,A,x_0)$ goes to the restriction $f|_{(I^{n-1}\times\qty{1}, \p I^{n-1}\times\qty{1})}$ under $\p_{\ast}$. 

\subsection{Cofiber sequence} \label{dis:cofib}
\noindent For any continuous map $f:X \to Y$ we know it can be decomposed as a cofibration and a homotopy equivalence. Cocide $Mf$ be the mapping cone over $Y$. Let, $j :X \to Mf$ be the inclusion $x \mapsto (x,1)$ and $r: Mf \to Y$ defined by $y \mapsto y$ and $(x,s) \mapsto f(x)$. Clearly $j$ is a cofibration and $r$ be a homotopy equivalence. For a based map $f:X \to Y$ we define homotopy cofibre $Cf$ to be $$Cf = Y\cup_f CX = Mf/j(X)$$ Let $\pi : Cf \to Cf/Y$ be the quotient map.

\[\begin{tikzcd}
	X & Y & Cf & {\Sigma X} & {\Sigma Y} & {\Sigma Cf} & \cdots
	\arrow["f", from=1-1, to=1-2]
	\arrow["i", from=1-2, to=1-3]
	\arrow["\pi", from=1-3, to=1-4]
	\arrow["{-\Sigma f}", from=1-4, to=1-5]
	\arrow["{-\Sigma i}", from=1-5, to=1-6]
	\arrow[from=1-6, to=1-7]
\end{tikzcd}\] is called \textit{cofiber sequence.} here, $-\Sigma f$ is the map that sends $(x \wedge t)$ to $f(x)\wedge (1-t)$. For the based spaces we have a definition of exactness. A sequence (let's say the above one) is said to be exact if the composition of two consecutive maps (excactness at $Y$: look at $i \circ f$) has image $\ast$(based point) iff it's pre image is only the based point.  The cofiber sequence turns out to be an excact sequence of spaces. Applying $[-,Z]$ functor we will get a long exact sequence of groups \[\begin{tikzcd}
	\cdots & {[\Sigma Cf,Z]} & {[\Sigma Y,Z]} & {[\Sigma X,Z]} & {[Cf,Z]} & {[Y,Z]} & {[X,Z]}
	\arrow[from=1-2, to=1-1]
	\arrow[from=1-3, to=1-2]
	\arrow[from=1-4, to=1-3]
	\arrow[from=1-5, to=1-4]
	\arrow[from=1-6, to=1-5]
	\arrow[from=1-7, to=1-6]
\end{tikzcd}\]

\noindent {\small \textcolor{black}{The above sequence is also known as Puppe sequence.}}
\subsection{Some Results}

\begin{itemize}
    \item $\pi_n(\s^n) =\Z$ for $n \geq 1$ and $\pi_i(\s^n) =0$ for $i<n$.
    \item For any fibration $E \to B$ with fibre being discrete/contractible, $\pi_n(E)=\pi_n(B)$ for all $n \geq 0$.
    \item  The `Hopf-fibration' $\s^1\hookrightarrow \s^3 \to \C P^1\simeq \s^2$ and thus $\pi_3(\s^2) =\pi_3(\s^3) \simeq \Z$.
    \item A triple $(X,A,B)$ is called \textbf{excisive triad} if $X = A^{\circ} \cup B^{\circ}$. For ordinary homology theory the inclusion of pairs $(A,A\cap B) \hookrightarrow (X,B)$ induces isomorphism in relative homology groups.
    \item For homotopy groups it's not the case. \textbf{Example} $X= \s^2\vee \s^2$ and $A = X_{+},B=X_{-}$ here $A \cap B \simeq \s^1 \vee \s^1$.
    \item \textbf{(Long exact sequence of triad)} If $(X,A,B)$ is a triad such that $B \subseteq A \subseteq X$ then we have the long exact sequence of relative homotopy groups $$\cdots \rightarrow \pi_n(A,B)\to \pi_n(X,B) \to \pi_n(X,A) \to \pi_{n-1}(A,B)\to \cdots$$
    \item[]  The above result will follow from chasing the following commutative diagram along \textcolor{red}{red arrow},\[\begin{tikzcd}
		{\pi_n(B)} & {\pi_n(A)} & {\pi_n(A,B)} & {\pi_{n-1}(B)} & {\pi_{n-1}(A)} & {} \\
		{\pi_n(B)} & {\pi_n(X)} & {\pi_n(X,B)} & {\pi_{n-1}(B)} & {\pi_{n-1}(X)} \\
		{\pi_n(A)} & {\pi_n(X)} & {\pi_n(X,A)} & {\pi_{n-1}(A)} & {\pi_{n-1}(X)} \\
		&&& {\pi_{n-1}(A,B)} & {\pi_{n-1}(X,B)}
		\arrow[from=1-1, to=1-2]
		\arrow[Rightarrow, no head, from=1-1, to=2-1]
		\arrow[from=1-2, to=1-3]
		\arrow[from=1-2, to=2-2]
		\arrow[from=1-3, to=1-4]
		\arrow[color={rgb,255:red,228;green,78;blue,115}, from=1-3, to=2-3]
		\arrow[from=1-4, to=1-5]
		\arrow[Rightarrow, no head, from=1-4, to=2-4]
		\arrow[from=1-5, to=2-5]
		\arrow[from=2-1, to=2-2]
		\arrow[from=2-1, to=3-1]
		\arrow[from=2-2, to=2-3]
		\arrow[Rightarrow, no head, from=2-2, to=3-2]
		\arrow[from=2-3, to=2-4]
		\arrow[color={rgb,255:red,228;green,78;blue,115}, from=2-3, to=3-3]
		\arrow[from=2-4, to=2-5]
		\arrow[from=2-4, to=3-4]
		\arrow[Rightarrow, no head, from=2-5, to=3-5]
		\arrow[from=3-1, to=3-2]
		\arrow[from=3-2, to=3-3]
		\arrow["{\p_{\ast}}", from=3-3, to=3-4]
		\arrow[color={rgb,255:red,228;green,78;blue,115}, from=3-3, to=4-4]
		\arrow[from=3-4, to=3-5]
		\arrow["{i_{\ast}}", from=3-4, to=4-4]
		\arrow[from=3-5, to=4-5]
		\arrow[color={rgb,255:red,228;green,78;blue,115}, from=4-4, to=4-5]
	\end{tikzcd}\]
\end{itemize}

\noindent We call a space \textbf{$n$-connected} if $\pi_n(X,x_0)\simeq 0$ and two space $X$ and $Y$ are \textbf{weakly equivalent} if $\pi_i(X) \simeq \pi_i(Y)$ for all $i \geq 0$. If $f:X \to Y$ is a map between two based spaces so that it induces isomorphism on every higher homotopy groups we call it an weak equivalence b/w the spaces. \ts{eg.} The sphere $\s^n$ is $(n-1)$ connected space. Weak equivalence may not be a homotopy equivalence. Consider $X = \qty{1/n, n \in \N}\cup\qty{0}$ and $Y$ is a countable discrete set. Then the natural map $f : X \to Y$ is weak equivalence but not homotopy equivalence. But this can be true if $X$ and $Y$ are CW complexes. 

\begin{theorem} \label{thm:whthd}
	(\textbf{Whitehead's Theorem}) If $f:X \to Y$ is a map between two connected CW complexes which is weak equivalence we can conclude $f$ is in-fact a homotopy equivalence. More generally if $A \hookrightarrow X$ is weak equivalence of CW complexes then $X$ deformation retract onto $A$.
\end{theorem}

\vspace*{0.2cm}

\noindent \textcolor{purple}{Whitehead's theorem doesn't say if two space $X$ and $Y$ have all homotopy groups same, then they are homotopic.} \ts{Example --} Let, $X = \R P^2$ and $Y = \s^2\times \R P^{\infty}$. These spaces are connected so $0$-th homotopy groups are same. $\pi_n(X) \simeq \pi_n(\s^2)$ for $n\geq 2$ and $\pi_n(Y) = \pi_n(\s^2) \times \pi_n(\R P^{\infty})\simeq \pi_n(\s^2)\times \pi_n(\s^{\infty}) \simeq \pi_n(\s^2)$. This is true for $n \geq 2$. Here we have used the fact there is a covering from the sphere to the projective space. Now for $g=1$ we can see $\s^2$ is simply connected and $\R P^2 \hookrightarrow \R P^{\infty}$ induces isomorphism on $\pi_1$. Thus, these space have same homotopy groups. \textbf{But}, the space $Y$ have non-trivial homology for infinitely many indexes unlike $X$. So $X \not\simeq_{\textbf{hTop}}Y$.

\vspace*{0.2cm}

We also can define \textbf{$n$-connectedness} of a pair $(X,A)$. A pair is said to be $n$-connected if $\pi_p(X,A)$ is trivial for $p\leq n$ and $\pi_0(X)\to \pi_0(A)$ is a surjection. Similarly, we can define \textbf{n-equivalence} of pairs. The map $f : (A,C)\to (X,B)$ (between connected spaces) is said to be $n$-equivalence if $f$ induces isomorphism in relative homotopy groups for indices $q <n$ and surjection for the index $q=n$. \ts{Eg:} For a CW complex $X$, the inclusion of $n$-th skeleton $X^n \hookrightarrow X$ is $n$-equivalence.

\subsection{CW approximations}

\noindent Here we list a few CW-approximation theorems we will state without proof. The proofs can be found at \cite[Page 76]{May}. 

\begin{itemize}
	\item[] (\textbf{Cellular Approximation}) Any map $f:(X,A) \to (Y,B)$ between pair of CW complexes is homotopic to a cellular map.
	\item[] (\textbf{Approximating a space by CW complex}) For any space $X$ there is a cellular complex $\Gamma X$ and a weak equivalence $\gamma : \Gamma X \to X$.  Such that given $f: X\to Y$ there is a map $\Gamma f : \Gamma X \to \Gamma Y$ so that the following diagram commutes, \[\begin{tikzcd}
		X & Y \\
		{\Gamma X} & {\Gamma Y}
		\arrow["f", from=1-1, to=1-2]
		\arrow["\gamma", from=2-1, to=1-1]
		\arrow["{{\Gamma f}}"', from=2-1, to=2-2]
		\arrow["\gamma"', from=2-2, to=1-2]
	\end{tikzcd}\]
	\item[] (\textbf{Approximating a pair by a pair of CW complex}) For any pair of spaces $(X, A)$ and any $C W$ approximation $\gamma$ : $\Gamma A \longrightarrow A$, there is a $C W$ approximation $\gamma: \Gamma X \longrightarrow X$ such that $\Gamma A$ is a subcomplex of $\Gamma X$ and $\gamma$ restricts to the given $\gamma$ on $\Gamma A$. If $f:(X, A) \longrightarrow(Y, B)$ is a map of pairs and $\gamma:(\Gamma Y, \Gamma B) \longrightarrow(Y, B)$ is another such $C W$ approximation of pairs, there is a map $\Gamma f:(\Gamma X, \Gamma A) \longrightarrow(\Gamma Y, \Gamma B)$, unique up to homotopy, such that the following diagram of pairs is homotopy commutative:
	
	\[\begin{tikzcd}
		{(X,A)} & {(Y,B)} \\
		{(\Gamma X, \Gamma A)} & {(\Gamma Y, \Gamma B)}
		\arrow["f", from=1-1, to=1-2]
		\arrow["\gamma", from=2-1, to=1-1]
		\arrow["{{\Gamma f}}"', from=2-1, to=2-2]
		\arrow["\gamma"', from=2-2, to=1-2]
	\end{tikzcd}\] If $(X, A)$ is n-connected, then $(\Gamma X, \Gamma A)$ can be chosen to have no relative $q$-cells for $q \leq n$.
\end{itemize}

\subsection{Eilenberg-MacLane Spaces}

Let $G$ be any group, and $n$ in $\mathbb{N}$. An \textit{Eilenberg-MacLane space} of type $(G, n)$, is a space $X$ of the homotopy type of a based CW-complex such that:
$$
\pi_{k}(X) \cong \begin{cases}G, & \text { if } k=n \\ 0, & \text { otherwise }\end{cases}
$$ One denotes such a space by $K(G, n)$. We now want to prove that the spaces $K(G, n)$ exist and are unique, up to homotopy, for every group $G$ and every integer $n \geq 0$. We will only show this statement when $G$ is an abelian group and when $n \geq 1$. The case $n=0$ is vacuous : one just takes the group $G$ endowed with its discrete topology. For $n=1$ and the group is not abelian $G$ can be represented as $G = \qty{\alpha_i : \beta_j}$ where $\alpha _i$ are generators and $\beta_j$ are relations. Now for each $\beta_j$ one disc should be attached to $\vee_{\alpha_i} \s^1$ according to the relation. This is how we can create a space with the homotopy group as required. Now we will use \textit{Homotopy killing} lemma to kill the higher homotopy groups by attaching $\geq 3$ cells. We know it don't affect the fundamental group \cite[chapter 1]{hatat}. Notice that when $n \geq 2$, the group $G$ must be abelian. Notice also that when $n \geq 1$, the spaces $K(G, n)$ are path-connected. More generally, the spaces $K(G, n)$ are $(n-1)$-connected.
\newcommand{\Lem}[1]{

      \vspace*{0.2cm}

      \noindent {${\S}$} \ts{Lemma --} \textit{#1}

      \vspace*{0.2cm}
    }

\Lem{\textsc{Homotopy Killing Lemma.} Let $X$ be any CW-complex and $n>0$. There exists a relative CW-complex $(X^{\prime}, X)$ with cells in dimension $(n + 1)$ only, such that $\pi_n(X^{\prime}) = 0$, and $\pi_k(X) \simeq \pi_k(X^{\prime})$ for $k < n$.}

\noindent \textit{Proof.} The proof is not very hard. But this idea will be very helpful. Let the generator of $\pi_n(X,x_0)$ are represented by $\qty{f_j:\s^n \to X: j \in \mathscr{J}}$. Here $\mathscr{J}$ is some index set. Consider the following pushout diagram \[\begin{tikzcd}
	{\bigvee_{j \in \mathscr{J}}\s^{n}} & X \\
	{\bigvee_{j \in \mathscr{J}}D^{n+1}} & {X'}
	\arrow["{\bigvee f_j}", from=1-1, to=1-2]
	\arrow[hook', from=1-1, to=2-1]
	\arrow["i", dashed, from=1-2, to=2-2]
	\arrow[dashed, from=2-1, to=2-2]
\end{tikzcd}\] Note that the map  $i$ is $n$-equivalence. For any generator $f_j$ the map $i \circ f_j : \s^n \to X'$ can be extended to a map $\hat{f}_j : D^{n+1} \to X^{\prime}$ (by the property of pushout). So the map $i \circ f_j$ is null-homotopic. In other words $\pi_n(i) : X \to X^{\prime}$ sends each generator $[f_j]$ to $[i \circ f_j]$ which is null homotopic. Thus the map $\pi_n(i)$ is trivial map since it is also surjective $\pi_n(X^{\prime})\simeq 0$. $\hfill \blacksquare$

\vspace*{0.2cm} \newcommand{\colim}{\text{colim}}

\noindent {\textbf{Existance}} of Eilenberg-MacLane spaces. We will show $k(\pi,n)$ exist for $n \geq 2$. For that we will consider \textit{Moore Spaces}. Briefly Moore-space $M(G,n)$ are the space which has integal simplicial homology $\simeq G$ for the index $n$ and trivial for other indices. It's not hars to show for abelian group $G$, Moore space always exist and infact by construction \cite[Example 2.40]{hatat} it is a CW complex. If $X = M(\pi,n)$ by the construction it don't have any cell of dimension $\leq (n-1)$ since $X_n \hookrightarrow X$ is $n$ -equivalence we can say $X$ is $(n-1)$ connected. By \ref{thm:hurw} we can say $\tilde{H}_n(X) \simeq \pi_n(X) \simeq \pi$ (as $n\geq 2$). We can construct a space $F_1X$ from $X$ by attaching $(n+2)$-cells so that $\pi_{n+1}(F_1X)$ is trivial. Iterate the process and by taking colimit $\colim_{j} F_jX$ we will get a space $\tilde{X}$ such that $\pi_k(\tilde{X})$ is trivial for $k>n$.

\vspace*{0.2cm}

\noindent \texttt{Another way -} There is a beutiful way to construct Eilenberg-MacLane spaces using`\textit{Infinite Symmetric Products}'. Given any based topological space $X$ we can construct a monoid $SP(X)$ in the following way: consider the action of $S_n$ on $X^n$ given by $\sigma . (x_1,\cdots,x_n) = (x_{\sigma(1)},\cdots, x_{\sigma(n)})$. Denote the orbit space of this action by $SP^n(X):= X^n /S_n$ (It can be shown it is functorial construction). Now there is a natual inclusion of $SP^n(X)\hookrightarrow SP^{n+1}(X)$ by $[x_1,\cdots,x_n] \to [x_1,\cdots,x_n,\ast]$ where $\ast$ is the based point of $X$. Define $$SP(X):= \colim \qty(\cdots SP^n(X)\hookrightarrow SP^{n+1}(X)\to \cdots)$$ Note that $SP^n(X)$ is a quotient space of $X^n$ thus it have a induced topology on it. We give $SP(X)$ the colimit topology i.e any subset $U \subseteq SP(X)$ is open iff $U \cap SP^n(X)$ is open for all $n$. By construction we can view $SP(X)$ as a topological space as well as a monoid (the product is the natural one with identity being $(\ast,\ast,\cdots)$). Now for any CW-complex $X$ we can give $SP(X)$ a CW-complex structure, at first we give CW complex structure to $X^n$ and then it induce a CW structure on $SP^n(X)$ and the colimit topology will help us to  get the CW-structure on $SP(X)$. \ts{Eg.} $SP (\s^2)= \C P^{\infty}$. Note that we can view $\C P^n$ as the equivalence class of polynomials $f(x)=a_nx^n + \cdots + a_0 \in \C[X]$ such that $f \simeq g \iff f = \lambda g$ for some complex scaler $\lambda$. We can view $\s^2$ as extended complex plane. So, $$SP^n(\s^{2})\to \C P^n : [a_1,\cdots ,a_n] \to [(x+a_1)(x+a_2)\cdots(x+a_n)]$$ is an bijection and by closed map lemma it is a \textit{homeomorphism}. Thus taking colimit will give us the result. \ts{Interesting Fact.} For any CW-complex $X$ there is a natural map $p:SP(X) \to \Omega(SP(\Sigma X))$ given by $[x_1,\cdots,x_{n=1}]\mapsto (t \mapsto [(x_1,t),\cdots,(x_n,t)])$. It turns out to be an weak equivalance. Thus we can conclude,$$\pi_n(SP(X)) = [S^n,SP(X)]_{\ast} = [S^n,\Omega SP(\Sigma X)]_{\ast} = [\s^{n+1}, SP(\Sigma X)]= \pi_{n+1}(SP(\Sigma X))$$ This helps us to define a \textit{homology theorey} for CW-complexes. Define $h_n(X,\ast):= \pi_n(SP(X),[\ast,\ast,\cdots])$. We can show, $\tilde{h}_n = h_n(-,\ast)$ satisfy three axioms \textit{Suspension axiom, Existance of Long exaact sequence of pairs, Additive axioms}. On the category of CW-complexes. If any homology theory staisfy these three axioms they are equivalent to \textit{the ordinary homology theory for CW-complexes}. Since we calso know any ordinary homology theoies are same we can say $\tilde{h}_n$ is infact equivalent to the Cellular homology theorey. Thus if we construct $M(\pi,n)$ (which is a $CW$-complex) we can say $SP (M(\pi,n))$ is $K(\pi,n)$. This is another way to see the Existance of Eilenberg MacLane spaces. 

\vspace*{0.2cm}

\noindent \textcolor{purple}{\textsc{Remark.}} Also $SP(X)$ can be thought as the commuatative version of james reduced product space. James product gives rise to a very Interesting monoid which have a nice cohomology ring structure and it also helps to give us EHP sequence \cite[section 2]{mynote1}. Now we will show if we restrict the definition of Eilenberg-MacLane spaces to only the CW-complexes we can prove it's unique upto homotopy equivalence. Till now we have worked in the category $\textbf{Top}_{\ast}$ (topological spaces with a base point). The discussion in section 1, the approximation theorems indicates it is enough if we deal the homotopy theory in the category of CW-complexes. From now onward we will cosider $\pi_n,n \geq 1$ to be a functor from the category of pointed CW-complexes \newcommand{\cw}{\ts{CW}_{\ast}}$\cw$ to \ts{Groups}.

\begin{theorem} 
	The Eilenberg MacLane spaces $K(\pi,n)$ are unique upto homotopy equivalence, where $\pi$ is abelian and $n\geq 1$.
\end{theorem}

\noindent For the proof we propose the following proposition/lemma. 

\begin{lemma}{}{}
	\hspace*{0.2cm} If $Y$ is a space such that $\pi_k(Y)$ is trivial for $k>n$ and $X$ be a CW-coplex with a subcomplex $X_{n+1}\subseteq A$. Then Any map $f: A \to Y$ can be extended to a map $\tilde{f}:X \to Y$.
\end{lemma}

\noindent \textit{Proof.} We can extend the map $f$ cell by cell. Consider $X_0 = A \cup e^{n+1}$ then for the attaching map $\varphi : \s^n \to A$, $f \circ \varphi$ can be extended to a map $\tilde{\varphi}: X_0 \to Y$ as it is null-homotopic. It is always the case for $k>n$. Soo we can Indeed extend the map to whole $X$. \ts{Note-} this is a general idea in \textit{obstruction theory}, that if a mapp can be extend to the whole space then it might represent something null-homotopic in the homotopy groups. Infact the subject obstruction theory is the study of possibilities od extending a map from a subcomplex to the wholw space. 

\begin{lemma}{}{}
	\hspace*{0.2cm} Let $Y= k(\pi,n)$ where $\pi$ is an abelian group and $n\geq 1$ and $X$ is a $(n-1)$-connected CW-complex. We have a natural map $\Phi : [X,Y]_{\ast} \to \text{Hom}_{\Z}(\pi_n(X),\pi)$ given by $$\Phi : [f]_{\ast}\mapsto \pi_n(f)$$ is bijection.
\end{lemma}

\noindent \textit{Proof.} As we have done previously, we will work with $(n+1)$-th skeleton of $X$ only. Since $X$ is $(n-1)$-conneceted we can assume $X_n$ is wedge of spheres $X _n = \bigvee_{J}\s^n$ and $X_{n+1}$ is given  by the following pushout \[\begin{tikzcd}
	{\bigvee_{j \in \mathscr{J}}\s^{n}} & {X_{n}} \\
	{\bigvee_{j \in \mathscr{J}}D^{n+1}} & {X_{n+1}}
	\arrow["\sum \varphi_{\alpha}", from=1-1, to=1-2]
	\arrow[hook', from=1-1, to=2-1]
	\arrow["i", dashed, from=1-2, to=2-2]
	\arrow[dashed, from=2-1, to=2-2]
\end{tikzcd}\] If $[f]_{\ast}$ and $[g]_{\ast}$ are two distinct equivalence class in $[X,Y]_{\ast}$ so that $\pi_n(f)=\pi_n(g)$, i.e $[f\circ h] = [g \circ h]$ for any map $h : \s^n \to X$. By surjectivity of $\pi_n(i)$ we can say there is a map $h': \s^n \to X_n$ so that $\pi_n(i)([h']) = [h]$. In particular $\pi_n(f\circ i)(h')= \pi_n(g \circ i)(h')$. If we consider the generators of $X_n$ as in \ref*{thm:hurw} they will have same image under $\pi_n(f \circ i) = \pi_n(g \circ i)$ where $h'$ represent generators on $\pi_n(X_n)$ given by the inclusion of the spheres in $X_n$. Thus $[f\circ i]_{\ast}=[g\circ i]_{\ast}$. Let, $H : X_n \times I \to Y$ be the homotopy b/w $f \circ i$ and $g\circ i$, it will help us to get a continuous map $\hat{H}: X_n\times I \cup X \times \p I:Y$ where $\hat{H}(X \times \p I) = f \cup g$. Now note that $(n+1)$ skeleton of $X\times I$ is $X_n\times I \cup X \times \p I$. So we can extend the homotopy to get  a homotopy $\tilde{H}: X \times I \to Y$ b/w $f$ and $g$. Thus $[f]_{\ast}=[g]_{\ast}$. It proves $\Phi$ is \textit{Injective}.

\vspace*{0.2cm}

\noindent Let $h: \pi_{n}(X) \rightarrow \pi_{n}(Y)$ be a group homomorphism. Let $X_{n}=\bigvee_{j \in J} S_{j}^{n}$. The group $\pi_{n}\left(X_{n}\right)$ is generated by the homotopy classes of the inclusions $\iota_{j}: S_{j}^{n} \hookrightarrow$ $\bigvee_{j \in \mathscr{J}} S_{j}^{n}$, 
$$
\pi_{n}\left(X_{n}\right) \xrightarrow{\pi_n(i)} \pi_{n}(X) \xrightarrow{h} \pi_{n}(Y),
$$ we define $f_{j}: S_{j}^{n} \rightarrow Y$ as a representative of the image of $\left[\iota_{j}\right]_{*}$, i.e. : $h\left(\pi_n(i)(\left(\left[\iota_{j}\right]_{*}\right))=\left[f_{j}\right]_{*}\right.$, for each $j$ in $\mathcal{J}$. The maps $\left\{f_{j}\right\}$ determine a map $f_{n}: X_{n} \rightarrow Y$ where $f_{n} \circ \iota_{j}=f_{j}$.  For each $\beta$ in $\mathscr{J}$, the map $i \circ \varphi_{\beta}$ is nullhomotopic. Hence $\pi_n\left(f_{n}\right)\left(\left[\varphi_{\beta}\right]_{*}\right)=h\left(\pi_n(i)(\left(\left[\varphi_{\beta}\right]_{*}\right))\right)=0$. Hence $f_{n} \circ \varphi_{\beta}$ is nullhomotopic for each $\beta$. Therefore $f_{n}$ extends to a map $f: X \rightarrow Y$, by the previous proposition. From $h \circ \pi_n(i)=\pi_n\left(f_{n}\right)=\pi_n(f)\circ \pi_n(i)$, since $i_{*}$ is surjective, we obtain that $\pi_n(f)=h$, i.e. : $\Phi\left([f]_{*}\right)=h$. Thus the function $\Phi$ is \textit{surjective}. $\hfill \blacksquare$

\vspace*{0.2cm}

\noindent Let $X$ and $Y$ be Eilenberg-MacLane space of type $(\pi, n)$. This means that there are isomorphisms $\theta: \pi_{n}(X) \rightarrow \pi$ and $\rho: \pi_{n}(Y) \rightarrow \pi$. From the previous theorem, the composite $\rho^{-1} \circ \theta: \pi_{n}(X) \rightarrow \pi_{n}(Y)$ is induced by a unique homotopy class of $X \rightarrow Y$ which is therefore a weak equivalence. Since $X$ and $Y$ are CW-complexes, the Whitehead Theorem implies that $X$ and $Y$ are homotopy equivalent. 

$\hfill \text{\small \textcolor{gray}{(end of the theorem)}}$

\vspace*{0.2cm}

\noindent \ts{Result (Milnor).} If $X$ is a based CW-complex then $\Omega X$ is also a based CW-complex. From here we can conclude $$\Omega K(\pi,n) \simeq_{\textbf{hTop}_{\ast}} K(\pi,n-1)$$not only that we can take adjuction to get $\Sigma K(\pi,n) \simeq_{\textbf{hTop}_{\ast}} K(\pi,n+1)$. 

\subsection{Homotopy Excision Theorem}
Theorem
In the previous section we have seen excision doesn’t hold for homotopy groups (unlike homology/-
cohomology groups). Thus it is difficult to compute the higher homotopy groups in this case neither
8
In the previous section we have seen excision doesn't hold for homotopy groups (unlike homology/cohomology groups). Thus it is difficult to compute the higher homotopy groups in this case neither we have Van-Kampen type of theorem. Homotopy excision theorem is the closest we can get in terms of excision for homotopy groups. 

\begin{theorem}
	(\textbf{Homotopy Excision/Blakers-Massey Theorem}) Suppose $(X,A,B)$ is excisive triad with $C= A \cap B$ such that $(A,C)$ is $n$-connected and $(B,C)$ is $m$-connected then the inclusion $i:(A,C) \hookrightarrow (X,B)$ induces isomorphism on relative homotopy groups $$i_{\ast}: \pi_q(A,C)\xrightarrow{\simeq} \pi_q(X,B)$$ for $q < m+n$ and it's surjection on relative homotopy groups for $q=m+n$ (In other words it is an $m+n$-equivalence)
\end{theorem}

\noindent\ts{Reduction 1 -} \textcolor{purple}{Enough to prove the statement for the triple $(X,A,B)$ where we get $A$ by attaching cells of dimension $>n$ to $C$ and we get $B$ by attaching cells of dimension $>m$ and $X = A\cup_{C}B$.} 

\begin{itemize}
	\item[] We can construct a pair $(\Gamma A, \Gamma C)$ such that $\gamma : \Gamma A \to A$ is weak equivalence and $\Gamma C$ is subcomplex of $\Gamma A$ such that $\gamma |_{\Gamma C}: \Gamma C \to C$ is also an weak equivalence. We can construct a CW complex $A_0$ from $\Gamma C$ by attaching cells so that the space $A_0$ and $\Gamma A$ have same homotopy groups for indices $>n$. We can do this by attaching cells of dimension $>n$ as $\pi_i(A) = \pi_i(\Gamma A) \simeq \pi_i(C)$ for $i \leq n$. By theorem \ref{thm:whthd}, we can say we can say $\Gamma A$ and $A_0$ are homotopic spaces. \[\begin{tikzcd}
		{\Gamma C} & {A_0} \\
		{\Gamma A}
		\arrow[hook, from=1-1, to=1-2]
		\arrow[hook', from=1-1, to=2-1]
		\arrow["{\simeq_{\textbf{hTop}}}", from=1-2, to=2-1]
	\end{tikzcd}\]
	The cells we have attached to $\Gamma C$ we will attach them to $C$ accordingly (with pre-composing with $\gamma |_{\Gamma C}$). Since everything we are doing upto weak equivalence it will be enough to deal with the reduction. 
\end{itemize}

\noindent\ts{Reduction 2 -} \textcolor{purple}{Enough to prove the statement for the triple $(X,A,B)$ where we get $A$ by attaching only one cells of dimension $>n$ to $C$ and we get $B$ by attaching only one cells of dimension $>m$ and $X = A\cup_{C}B$.} 

\begin{itemize}
	\item[] Assume $A$ and $B$ are constructed by attaching cells as we have described in the first reduction (we are dealing with CW approximations only but renaming them with the initial characters only). Consider $C \subset A' \subset A$ such that $A$ is obtained from $A'$ by attaching one cells. Now as a pair $(A,C)$ has one more cell then $(A',C)$. Consider $X'= A'\cup_C B$. If excision holds for $(X',A',B)$ and $(X,X',B)$ then from the following commutative diagram (using five lemma) we get, excision holds for $(X,A,B)$ too. \[\begin{tikzcd}
		{\pi_{k+1}(A,A')} & {\pi_k(A',C)} & {\pi_k(A,C)} & {\pi_k(A,A')} & {\pi_{k-1}(A',C)} \\
		{\pi_{k+1}(X,X')} & {\pi_{k}(X',B)} & {\pi_k(X,B)} & {\pi_k(X,X')} & {\pi_{k-1}(X',B)}
		\arrow[from=1-1, to=1-2]
		\arrow[Rightarrow, no head, from=1-1, to=2-1]
		\arrow[from=1-2, to=1-3]
		\arrow[Rightarrow, no head, from=1-2, to=2-2]
		\arrow[from=1-3, to=1-4]
		\arrow[from=1-3, to=2-3]
		\arrow[from=1-4, to=1-5]
		\arrow[Rightarrow, no head, from=1-4, to=2-4]
		\arrow[Rightarrow, no head, from=1-5, to=2-5]
		\arrow[from=2-1, to=2-2]
		\arrow[from=2-2, to=2-3]
		\arrow[from=2-3, to=2-4]
		\arrow[from=2-4, to=2-5]
	\end{tikzcd}\]
\end{itemize}

\noindent \textbf{Proof for the reduced case.} \noindent Let, $A = C \cup e, B = C \cup e'$, choose $x \in e^{\circ}, y\in e'^{\circ}$. note that $\pi_i(A,C) \simeq \pi_{i}(X-y,X-\qty{x,y})$ and $\pi_i(A,C) \simeq \pi_{i}(X-y,X-\qty{x,y})$. These isomorphisms follow from the observation that \(X \backslash\{x\}\) is homotopy equivalent to \(B\) by retracting \(e \backslash\{x\}\) to its boundary, with similar retractions yielding \(X \backslash\{y\} \simeq A\) and \(X \backslash\{x, y\} \simeq C\).  

\vspace{0.2cm}

\noindent We first discuss surjectivity. Consider a representative of \(\pi_i(X, B)\), that is, a map \(f:\left(I^i, \partial I^i\right) \rightarrow (X, B)\) that takes \(J^{n-1}\) to the basepoint \(* \in C\). In other words, \(f\) maps the top face of \(I^i\) into \(B\) and the rest of the boundary to \(*\). By the diagram above, it suffices to prove that \(f\) is homotopic to a map \(f^{\prime}\) via a homotopy \(h\), such that  
(i) the image of \(f^{\prime}\) is in \(X \backslash\{y\}\),  
(ii) for every \(t \in I\), the restriction of \(h_t\) to the top face of \(I^i\) avoids \(x\),  
(iii) for every \(t \in I\), \(h_t\) maps \(J^{i-1}\) to \(*\).  If such a map \(f^{\prime}\) and homotopy \(h\) exist, then every representative in \(\pi_i(X, B) \cong \pi_i(X, X \backslash\{x\})\) is homotopic to some representative in \(\pi_i(X \backslash\{y\}, X \backslash\{x, y\})\), establishing that \(\pi_i(A, C) \rightarrow \pi_i(X, B)\) is surjective for \(i \leq n+m\). For the proof that such \(f^{\prime}\) and \(h\) exist, we refer the reader to \cite{May}.  

\vspace{0.2cm}

\noindent Injectivity of \(\pi_i(A, C) \rightarrow \pi_i(X, B)\) for \(i < n+m\) follows by an analogous argument. Suppose two representatives \(g, g^{\prime}\) of \(\pi_i(A, C)\) satisfy \([g] = [g^{\prime}] \in \pi_i(X, B)\) via a homotopy \(H: I^i \times I \rightarrow X\). Replacing \(f\) in the previous argument with \(H\), we claim that there exists a new map \(H^{\prime}\), homotopic to \(H\) via a homotopy \(G\), such that \(H^{\prime}\) avoids \(y\), the restriction of \(G_t\) to the top face of \(I^{i+1}\) avoids \(x\), and \(G_t\) maps \(J^i\) to \(*\). This establishes a homotopy from \(f\) to \(g\) in \(X \backslash\{y\}\), relative to \(X \backslash\{x, y\}\), implying \([g] = [g^{\prime}] \in \pi_i(X \backslash\{y\}, X \backslash\{x, y\})\). Since this holds for \(i+1 \leq n+m\) (where the domain of \(H\) is the \((i+1)\)-cube and the domain of \(f\) is the \(i\)-cube), we conclude injectivity for \(i < n+m\). $\hfill \blacksquare$

\subsection{Freudenthal Suspension Theorem}

If $f$ is a based map $f : \s^k \to X$ the suspension $\Sigma f : \Sigma \s^k \to \Sigma X$ given by $(x \wedge t) \mapsto f(x)\wedge t$. As we have already said suspension (reduced suspension) is a functor from $\textbf{Top}_{\ast}$ to itself. From the above discussion we see $\Sigma$  gives us a map $\pi_k(X) \to \pi_{k+1}(\Sigma X)$. Infact we can view $\Sigma : \pi_{k} \Rightarrow \pi_{k+1}$ as a natural transformation as the following diagram commutes for any $f :  X\to Y$
 \[\begin{tikzcd}
	{\pi_k(X)} & {\pi_{k+1}(\Sigma X)} \\
	{\pi_k(Y)} & {\pi_{k+1}(\Sigma Y)}
	\arrow["\Sigma", from=1-1, to=1-2]
	\arrow["{f_{\ast}}"', from=1-1, to=2-1]
	\arrow["{(\Sigma f)_{\ast}}", from=1-2, to=2-2]
	\arrow["\Sigma"', from=2-1, to=2-2]
\end{tikzcd}\]

\noindent Since we have excision kind of tools for computing homotopy groups we will establish the Freudenthal suspension theorem it will help us to get idea about stable homotopy theory.

\begin{theorem}\label{thm:freu}
	(\textbf{Freudenthal Suspension Theorem}) If $X$ is a $(n-1)$ conneceted CW complex, the map $\Sigma : \pi_{k}(X) \to \pi_{k+1}(\Sigma X)$ is an isomorphism for $k\leq 2n-2$ and surjective for $n = 2n-1$. 
\end{theorem}

\textit{Proof.} Let, $X$ be the based space with based point $x_0$ then we can view $\Sigma X$ as the pushout of the following diagram,
\[\begin{tikzcd}
	{X \times\p I \cup \qty{x_0} \times I} & {X\times I} \\
	{\qty{x_0}} & {\Sigma X}
	\arrow[hook, from=1-1, to=1-2]
	\arrow[from=1-1, to=2-1]
	\arrow[from=1-2, to=2-2]
	\arrow[from=2-1, to=2-2]
\end{tikzcd}\]

\noindent Consider the open cover of $X$, $A = X\times (0,1]/X \times \qty{1} \cup \qty{x_0}\times \left(0,1\right]$ and $B=X\times [0,1)/X \times \qty{0} \cup \qty{x_0}\times \left[0,1\right)$. We can see that $A$ and $B$ are open in $\Sigma X$, and there are the based homotopy equivalences $A \simeq_* C X, B \simeq_* C^{\prime} X, A \cap B \simeq_* X$. Where $CX$ and $C'X$ are reduced cone on $X$ defined by $C'X = X\times [0,1]/X \times \qty{0} \cup \qty{x_0}\times \left[0,1\right]$ and $CX = X\times [0,1]/X \times \qty{1} \cup \qty{x_0}\times \left[0,1\right]$. Indeed, the homotopy,
$$ 
\begin{aligned}
H: C X \times I & \longrightarrow C X \\
([x, t], s) & \longmapsto[x, s+(1-s) t],
\end{aligned}
$$
gives a based homotopy equivalence $C X \simeq_* x_0$. With the same argument, we have :
$$
A \simeq_* x_0, \text { and } B \simeq_* x_0 \simeq_* C^{\prime} X
$$

\noindent Hence, the triad $(\Sigma X ; A, B)$ is excisive. Moreover, $A$ and $B$ are contractible spaces, so $(A, X)$ and $(B, X)$ are $(n-1)$-connected, by the long exact sequence of the pairs, whence we can apply the excision homotopy Theorem. The inclusion $(B, X) \hookrightarrow(\Sigma X, A)$ is a $(2 n-2)$-equivalence, and thus, the inclusion $i:\left(C^{\prime} X, X\right) \hookrightarrow(\Sigma X, C X)$ is a $(2 n-2)$-equivalence. To end the proof, we need to know the relation between the inclusion $i$ and the suspension homomorphism $\Sigma$. Consider an element $[f]_* \in \pi_k(X)=\left[\left(I^k, \partial I^k\right),(X, *)\right]_*$. Let us name $q: X \times I \rightarrow C^{\prime} X \cong X \times I /(X \times\{0\} \cup\{*\} \times I)$ the quotient map induced by the definition of $C^{\prime} X$ as a pushout. Define $g$ to be the composite :
$$
I^{k+1} \xrightarrow{f \times \text { id }} X \times I \xrightarrow{q} C^{\prime} X .
$$

It is easy to see that $g\left(\partial I^{k+1}\right) \subseteq X$, and $g\left(J^k\right)=\{*\}$. Indeed, we have $g\left(\partial I^k \times I\right)=[(*, I)] \subseteq X$, $g\left(I^k \times\{0\}\right)=* \in X$ and $g\left(I^k \times\{1\}\right) \subseteq X$. Hence $g\left(\partial I^{k+1}\right) \subseteq X$. It is similar to prove that $g\left(J^k\right)=\{*\}$. Therefore $[g]_* \in \pi_{k+1}\left(C^{\prime} X, X\right)$. Moreover, it is clear $g{|_{I^k \times\{1\}}}=f$. Hence $\partial\left([g]_*\right)=[f]_*$, where $\partial$ is the boundary map of the long exact sequence of the pair $\left(C^{\prime} X, X\right)$. We get : $\rho \circ g=\Sigma f$, where the map $\rho: C^{\prime} X \rightarrow \Sigma X$ can be viewed as a quotient map, through the homeomorphism $\Sigma X \cong C^{\prime} X /(X \times\{1\})$. Thus the following diagram commutes :
\[\begin{tikzcd}
	{\pi_{k+1}(C'X,X)} & {\pi_k(X)} \\
	{\pi_{k+1}(\Sigma X, CX)} & {\pi_{k+1}(\Sigma X)}
	\arrow["\p", from=1-1, to=1-2]
	\arrow["{i_{\ast}}"', from=1-1, to=2-1]
	\arrow["\rho", from=1-1, to=2-2]
	\arrow["\Sigma", from=1-2, to=2-2]
	\arrow["\simeq", from=2-2, to=2-1]
\end{tikzcd}\]

\noindent Here $\p$ (from LES of pair $(C'X,X)$) and $i_{\ast}$(from the excision theorem) are also isomorphism. Thus $\Sigma$ is also an isomorphism for $k < 2n-1$ and surjection for $k =2n-1$ as $i_{\ast}$ is a surjection.

\begin{theorem} \label{thm:1}
	Let $f: X \longrightarrow Y$ be an $(n-1)$-equivalence between $(n-2)$-connected spaces, where $n \geq 2$; thus $\pi_{n-1}(f)$ is an epimorphism. Then the quotient map $\pi:(M f, X) \longrightarrow(C f, *)$ is a $(2 n-2)$-equivalence. In particular, $C f$ is $(n-1)$ connected. If $X$ and $Y$ are $(n-1)$-connected, then $\pi:(M f, X) \longrightarrow(C f, *)$ is a $(2 n-1)$-equivalence.
\end{theorem}

\noindent \textit{Proof.} We are writing $C f$ for the unreduced cofiber $M f / X$. We have the excisive $\operatorname{triad}(C f ; A, B)$, where $$
A=Y \cup(X \times[0,2 / 3]) \quad \text { and } \quad B=(X \times[1 / 3,1]) /(X \times\{1\}) .
$$ Thus $C \equiv A \cap B=X \times[1 / 3,2 / 3]$. It is easy to check that $\pi$ is homotopic to a composite
$$
(M f, X) \stackrel{\simeq}{\simeq}(A, C) \longrightarrow(C f, B) \stackrel{\simeq}{\longrightarrow}(C f, *),
$$
the first and last arrows of which are homotopy equivalences of pairs. The hypothesis on $f$ and the long exact sequence of the pair $(M f, X)$ imply that $(M f, X)$ and therefore also $(A, C)$ are $(n-1)$-connected. In view of the connecting isomorphism $\partial: \pi_{q+1}(C X, X) \longrightarrow \pi_q(X)$ and the evident homotopy equivalence of pairs $(B, C) \simeq(C X, X),(B, C)$ is also $(n-1)$-connected, and it is $n$-connected if $X$ is $(n-1)$-connected. The homotopy excision theorem gives the conclusions.

\subsection{Hurewicz Theorem}

There is a special relation between the homotopy groups of a space and the ordinary homology groups of that space with integral coefficients. Infact we can naturally produce a map $h: \pi_n(X) \to H_n(X,x_0)$ (here we are dealing with the based space $X$). We know the relative homology groups of $(\s^n,e)$ with integral coefficients is isomorphic to $\Z$. We can assume $i_n$ to be the generator of $H_n(\s^n,e)$. Then $h([f]) = H_n(f)(i_n)$ is a well defined map from homotopy group to relative homology group as homology groups are homotopy invariant. It turns out to be a homomorphism b/w the groups and we call it \textit{Hurewicz homomorphism}. If $[f]$ and $[g]$ are two class of maps in $\pi_n(X)$ then $[f]+[g]=[f\vee g \circ c]$ where $c$ is the pinching map. From the following commutative diagram \[\begin{tikzcd}
	{\tilde{H}_n(\s^n)} & {\tilde{H}_n(\s^n \vee \s^n)} & {\tilde{H}_n(X)} \\
	& {\tilde{H}_n(\s^n) \oplus\tilde{H}_n(\s^n)}
	\arrow["{\tiny{\tilde{H}_n(c)}}", from=1-1, to=1-2]
	\arrow["\D"', from=1-1, to=2-2]
	\arrow["{\tiny{\tilde{H}_n(f\vee g)}}", from=1-2, to=1-3]
	\arrow["\simeq"', from=1-2, to=2-2]
	\arrow["{\tiny{\tilde{H}_n(f) + \tilde{H}_n(g)}}"', from=2-2, to=1-3]
\end{tikzcd}\]
we get $H_n(f\vee g \circ c)(i_n) = H_n(f \vee g)\circ H_n(c) =H_n(f)(i_n)+H_n(g)(i_n)$. Thus $h([f]+[g])= h(f)+h(g)$ and hence it is a group homomorphism. The homomorphism can be viewed as a natural functor $h : \pi_n \Rightarrow \tilde{H}_n$ for $n \geq 0$ and furthermore it's compatible with the suspension homomorphism i.e. the following diagram commutes,\[\begin{tikzcd}
	{\pi_n(X)} & {\tilde{H}_{n}(X)} \\
	{\pi_{n+1}(\Sigma X)} & {\pi_{n+1}(\Sigma X)}
	\arrow["h", from=1-1, to=1-2]
	\arrow["\Sigma"', from=1-1, to=2-1]
	\arrow["{\Sigma }", from=1-2, to=2-2]
	\arrow["h"', from=2-1, to=2-2]
\end{tikzcd}\]

\noindent \textcolor{purple}{\textsc{Remark}.} The Hurewicz homomorphism can be defined for any ordinary homology theories. Ordinary homology theories satisfy Eilenberg-steenrod axioms \cite[page 95]{May}. Any ordinary homology theories are equivalent. We will work with cellular homology for the rest part. 

\Lem{Consider the wedge of $n$-spheres $X^n=\bigvee_{j\in \mathscr{J}}\s^n$ where $\mathscr{J}$ is any index set, $i_j^n$ be the inclusion of $j$-th index sphere in the wedge. Then $\pi_1(X^1)$ is free group generated by $\qty{i^1_j}$ and for $n \geq 2$, $\pi_n(X^n)$ is free abelian group generated by $\qty{i_j^n}$.}

\noindent \textit{Proof.} The $n=1$ case follows from the Seifert Van Kampen theorem. Let us prove now the case $n \geq 2$. Let $\mathscr{J}$ be a finite set. Regard $\bigvee_{j \in \mathscr{J}} S^n$ as the $n$-skeleton of the product $\prod_{j \in \mathscr{J}} S^n$, where again the $n$-sphere $S^n$ is endowed with its usual CW-decomposition, and $\prod_{j \in \mathcal{J}} S^n$ has the CW-decomposition induced by the finite product of CW-complexes. Since $\prod_{j \in \mathscr{J}} S^n$ has cells only in dimensions a multiple of $n$, the pair $\left(\prod_{j \in \mathscr{J}} S^n, \bigvee_{j \in \mathscr{J}} S^n\right)$ is $(2 n-1)$ connected. The long exact sequence of this pair gives the isomorphism :
$$
\pi_n\left(\bigvee_{j \in \mathcal{J}} S^n\right) \cong \pi_n\left(\prod_{j \in \mathscr{J}} S^n\right) \cong \bigoplus_{j \in \mathcal{J}} \pi_n\left(S^n\right)
$$
induced by the inclusions $\left\{\iota_j^n\right\}_{j \in \mathscr{J}}$. The result follows. Let now $\mathscr{J}$ be any index set, let $\Theta_{\mathscr{J}}: \bigoplus_{j \in \mathscr{J}} \pi_n\left(S^n\right) \rightarrow \pi_n\left(\bigvee_{j \in \mathscr{J}}\right)$ be the homomorphism induced by the inclusions $\left\{\iota_j^n\right\}_{j \in \mathscr{J}}$. Just as the case $n=1$, one can reduce $\mathscr{J}$ to the case where it is finite to establish that $\Theta_{\mathscr{J}}$ is an isomorphism.  $\hfill \blacksquare$

\vspace*{0.2cm}

\noindent \textcolor{purple}{From the above lemma we conclude Hurewicz homomorphism $h: \pi_n(X^n)\to \tilde{H}_n(X^n)$ is isomorphism for $X^n, n\geq 2$ and it is the abelianization homomorphism for $n=1$. }Infact for any $(n-1)$ connected space the Hurewicz homomorphism is an isomorphism ($n>1$). This is the statement of Hurewicz isomorphism theorem.

\begin{theorem}\label{thm:hurw}
	Let $X$ be a $(n-1)$-connected based space, where $n \geq 1$ then the Hurewicz isomorphism,
	$$h: \pi_n(X)\to \tilde{H}_n(X)$$
	is the abelianization homomorphism if $n = 1$ and is an isomorphism if $n>1$. 
\end{theorem}

\noindent \textit{Proof.} (We will deal with $n \geq2$ at first) The weak equivalence induce isomorphism in the relative homology groups. By the CW approximation we can assume $X$ to be weak equivalent to the CW complex $\Gamma X$. It is enough to work with $\Gamma X$, it is also $(n-1)$ connected. By the whitehead approximation theorem we can assume $\Gamma X$ is homotopic to a CW complex $X'$ that do not have any cells of dimension $k$ (here $1\leq k \leq (n-1)$) and have one $0$-cell. Call this space $FX$. Since we are working with based spaces (i.e. $FX$ is bases space) the $n$-th skeleton of $FX$ is achived by the following pushout, 

\[\begin{tikzcd}
	{\bigvee_{j \in \mathscr{J}_n}\s^{n-1}} & {FX_0} \\
	{\bigvee_{j \in \mathscr{J}_{n}}D^n} & {FX_n}
	\arrow["g", from=1-1, to=1-2]
	\arrow[hook', from=1-1, to=2-1]
	\arrow[dashed, from=1-2, to=2-2]
	\arrow[dashed, from=2-1, to=2-2]
\end{tikzcd}\]

\noindent (Here $FX_k$ means $k$-th skeleton of $FX$) Thus $FX_n$ is nothing but wedge of spheres i.e $FX_n =\bigvee_{j\in\mathscr{J}_n}\s^n$. The $(n+1)$-the skeleton will also be constructed by similar kind of pushout. Note that, cone over $\bigvee_{j \in \mathscr{J}_{n+1}}\s^{n}$ is homeomorphic to wedge of disks and thus $FX_{n+1}$ is actually the mapping cone(reduced) over $f$. Where $f : \bigvee_{j \in \mathscr{J}_{n+1}}\s^{n} \to FX_{n}$ is the attaching map. The following diagram shall describe it clearly,

\[\begin{tikzcd}
	& {\bigvee_{j \in \mathscr{J}_{n+1}}\s^{n}} & {FX_n} \\
	{C \qty(\bigvee_{j \in \mathscr{J}_{n+1}}\s^{n})} & {\bigvee_{j \in \mathscr{J}_{n+1}}D^{n+1}} & {FX_{n+1}} & Cf
	\arrow["f", from=1-2, to=1-3]
	\arrow[hook', from=1-2, to=2-1]
	\arrow[hook', from=1-2, to=2-2]
	\arrow[dashed, from=1-3, to=2-3]
	\arrow["{\simeq_{\text{hTop*}}}"', from=2-1, to=2-2]
	\arrow[dashed, from=2-2, to=2-3]
	\arrow["{\simeq_{\text{hTop*}}}"', from=2-3, to=2-4]
\end{tikzcd}\]

\noindent As we know $FX_{n+1}\hookrightarrow FX$ is $n$-equivalence it will induce isomorphism if $n$-th homotopy, i.e. $\pi_n(FX_{n+1}) = \pi_n(FX)$ also from cellular homology theory \cite[page 137]{hatat} we know this inclusion will induce isomorphism on $n$-th reduced homology i.e. $\tilde{H}_n(FX_{n+1}) \simeq \tilde{H}_{n}(FX)$. Thus it is enough to prove the `Hurewicz isomorphism' for $FX_{n+1}$. Let us call the space $\bigvee_{j \in \mathscr{J}_{n+1}}\s^{n}:=T$. Thus we have a long exact sequence of homology groups \cite[page 128]{cloh}, $$\cdots \to \tilde{H}_n(T) \to \tilde{H}_n(FX_n) \to \tilde{H}_n(Cf)\to \tilde{H}_{n-1}(T)\to \cdots$$ Since $T$ is wedge of $n$-spheres $\tilde{H}_{n-1}(T)$ is trivial. In the following commutative diagram the bottom row is exact

\[\begin{tikzcd}
	{\pi_n(T)} & {\pi_{n}(FX_{n})} & {\pi_n(FX_{n+1})} & 0 \\
	{\tilde{H}_n(T)} & {\tilde{H}_n(FX_n)} & {\tilde{H}_n(\underbrace{FX_{n+1}}_{Cf})} & 0
	\arrow[from=1-1, to=1-2]
	\arrow["h"', from=1-1, to=2-1]
	\arrow[from=1-2, to=1-3]
	\arrow["h", from=1-2, to=2-2]
	\arrow[from=1-3, to=1-4]
	\arrow["h", from=1-3, to=2-3]
	\arrow[from=2-1, to=2-2]
	\arrow[from=2-2, to=2-3]
	\arrow[from=2-3, to=2-4]
\end{tikzcd}\]

\noindent Since $T$ and $FX_n$ are $(n-1)$ connected space and $f$ is $(n-1)$-equivalence, by \ref{thm:1} we must have $(Mf,T) \to (Cf,\ast)$ is $(2n-1)$-equivalence, so it induces isomorphism on $\pi_n$. Now using the LES for homotopy groups we get \[\begin{tikzcd}
	{\cdots } & {\pi_n(T)} & {\pi_{n}(Mf)} & {\pi_n(Mf,T)} & \cdots \\
	& {\pi_n(T)} & {\pi_n(FX_n)} & {\pi_n(Cf)}
	\arrow[from=1-1, to=1-2]
	\arrow[from=1-2, to=1-3]
	\arrow[Rightarrow, no head, from=1-2, to=2-2]
	\arrow[from=1-3, to=1-4]
	\arrow["\simeq", from=1-3, to=2-3]
	\arrow[from=1-4, to=1-5]
	\arrow["\simeq", from=1-4, to=2-4]
	\arrow[from=2-2, to=2-3]
	\arrow[from=2-3, to=2-4]
\end{tikzcd}\]
the top row will also be exact. By the previous lemma 1st and seconf Hurewicz homomorphism are isomorphism so we have proved the Hurewicz isomorphism for $FX_{n+1}\simeq Cf$. 

\subsection{Stability}

\noindent Let $X$ be a $(n-1)$ -connected space. We get the sequence following homotopy groups (by applying suspension) consecutively,

$$\pi_k(X)\xrightarrow{\Sigma} \pi_{k+1}(\Sigma X)\xrightarrow{\Sigma}\cdots \pi_{k+r}(\Sigma^r X)\xrightarrow{\Sigma}\cdots$$

\noindent Inductively we can show $\Sigma^rX$ is $(n+r-1)$-connected. Thus for larger $r$ the homotopy groups finally gets stabilized. We can define stable homotopy groups as follows,

\begin{definition}
	Let $X$ be a $(n-1)$-connected space. Let, $k \geq 0$ and the $k$-th \textbf{stable homotopy group} of $X$ is the colimit of $\pi_{k+r}(\Sigma^rX)$, $$\pi_k^S(X) := \underset{r}{\colim} \, \pi_{k+r}(\Sigma^r X)$$
\end{definition}

\noindent If $k<n-1$ we must have, $$\pi_k^S(X)= \pi_{k+n}(\Sigma^nX)$$ It is one of the interesting problem in algebraic topology to compute $\pi_k^s(\s^0)$. It arises in computation of different geometric things such as parallelizable structures on $\s^n$ for $n \geq 5$. In general the groups $\pi_{k+n}(\s^n)$ are called \textbf{stable} if $n>k+1$ and \textbf{unstable} if $n \leq k+1$.

\begin{theorem} \label{stab:fint}
	The stable homotopy groups of sphere are finite. If we define $\pi^S_k(\s^0)=: \pi^S_k$, this is finite. 
\end{theorem}

\noindent \textit{Proof.} The proof is technical but we will use a theorem by J.P.Serre, called \textit{Serre finiteness} \cite[Page 12]{mynote1}. This asserts the higher homotopy groups of $\s^n$ are finite except for the index $=n$ and the case $n=2k$ and index is $4k-1$. If $k\neq 1$ it's not hard to see that $\pi_k^S$ is finite. For $k=1$ we have a very nice result that, $$\pi_1^S\simeq \Z/2\Z$$ which is finite. We will prove that $\pi_4(\s^3) \simeq \Z/2\Z$ and using Freudenthal suspension theorem successively we will get, $$\pi_4(\s^3)=\pi_5(\s^4)=\pi_6(\s^5) \cdots = \pi_{n+1}(\s^n)$$ taking colimit will give us $\pi_1^S = \Z / 2\Z$ and thus the proof is complete. $\hfill \blacksquare$.

\section{Framed Cobordism-The Pontryagin Construction}

We know for any manifold $M$ and $N$ of same dimension, if we have a compactly supported map $f:M \to N$ then it induces a map in compactly supported cohomology $$f^{!}: H^n_c(N)\to H_c^n(M)$$ since the compactly supported top-degree cohomology are isomorphic to $\R$, $f^{!}$ is actually a linear map from$\R$ to $\R$. Thus image of this map is determined by $f^{!}(1)$. Interestingly it will turns out to be an integer. We call that integer, the degree of the map $f$. However, if we have a map which is not compactly supported we can't guarantee this. If the manifolds were closed and oriented then by Poincar\'e duality we can say their top cohomology is also isomorphic to $\R$ and similar idea will help us do define the degree of a map. The Pontryagin construction helps us to talk about the degree of maps $f:M \to \s^n$ for any compact boundaryless manifold. In the following definiton we will assume $\p M = \p N = \p N' = \emptyset$.

\begin{definition}
  A manifold $N$ is \textit{cobordant} to $N'$ within $M$ if the subset $N \times [0,\ep) \cup N' \times (1-\ep,1]$ of $M\times[0,1]$ canbe extended to a compact manifold $X \subset M \times [0,1]$ so that $$\p X = N \times \qty{0} \cup N' \times \qty{1}$$ and $X$ does not intersect $\p (M \times [0,1])$ except for $\p X$.
\end{definition}

\noindent If two sub-manifold $N,N'$ of $M$ are cobordant we will denote $N\sim_c N'$. It's  not hard to see it is an equivalence relation (the transitivity is shown in the following diagram)
\[\includegraphics[width = 9cm]{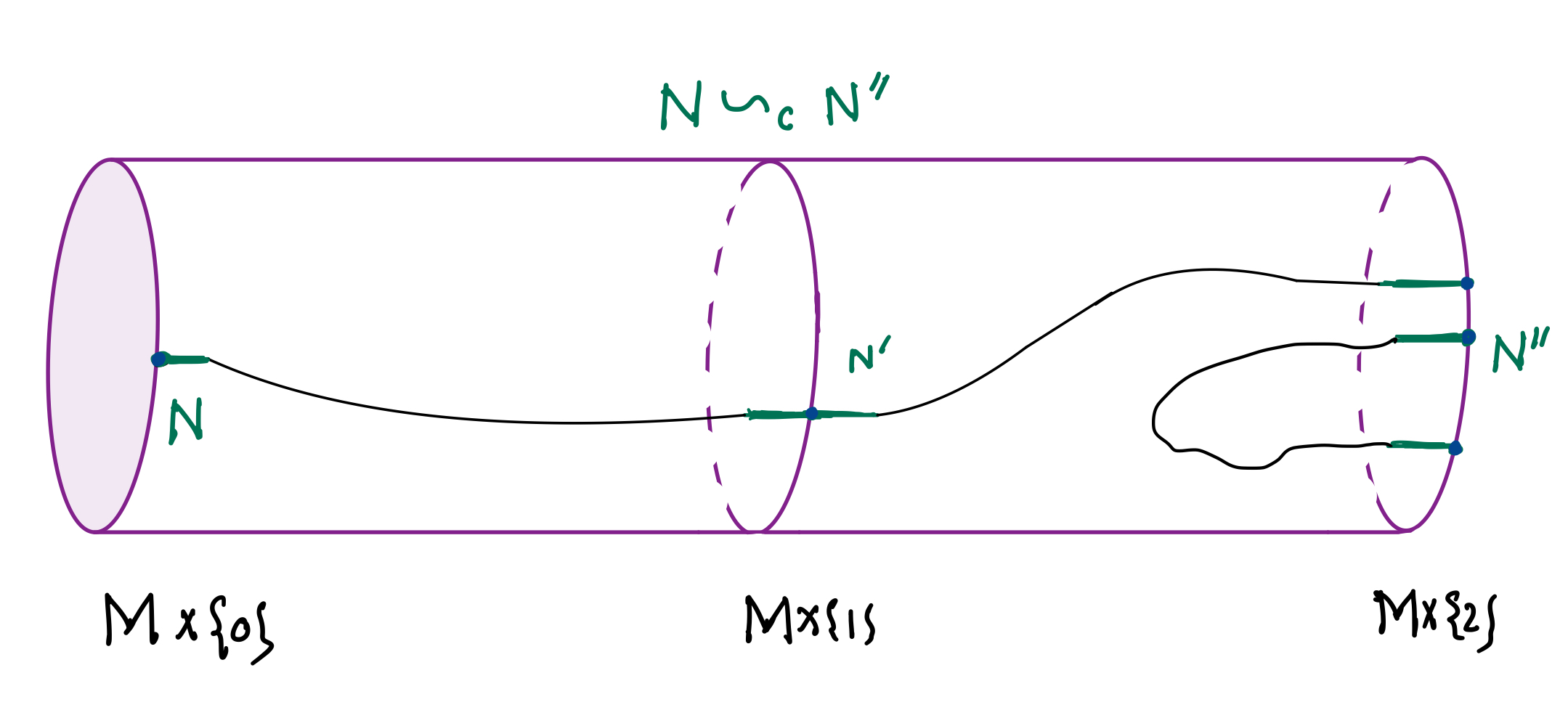}\]

\noindent Recall that \textit{framing} of a submanifold $N \subset M$ is a smooth function $\nu : N \to ((T_xN)^{\perp_{T_xM}})^{m-n}$ such that $\nu(x) = \qty(\nu^1(x),\cdots, \nu^{m-n}(x))$ is a basis of the orthogonal component of $T_xN$ inside $T_xM$, here $m-n$ is codimension of $N$ in $M$. The pair $(N,\nu)$ is called \textit{framed submanifold}. Two framed submanifold $(N,\nu)$ and $(N',\nu')$ are said to be \textit{framed cobordant} if therse exist a cobordism $X\subset M \times [0,1]$ between $N$ and $N'$ and a framing $u$ of $X$ such that (as shown in the picture)\begin{align*}
  u^i(x,t) = (\nu^i(x),0) & \text{ for } (x,t) \in N\times [0,\ep) \\
  u^i(x,t) = (\nu'^i(x),0) & \text{ for } (x,t) \in N \times (1-\ep,1]
\end{align*} It is also an equivalence relation. 

\[\includegraphics[width = 7cm]{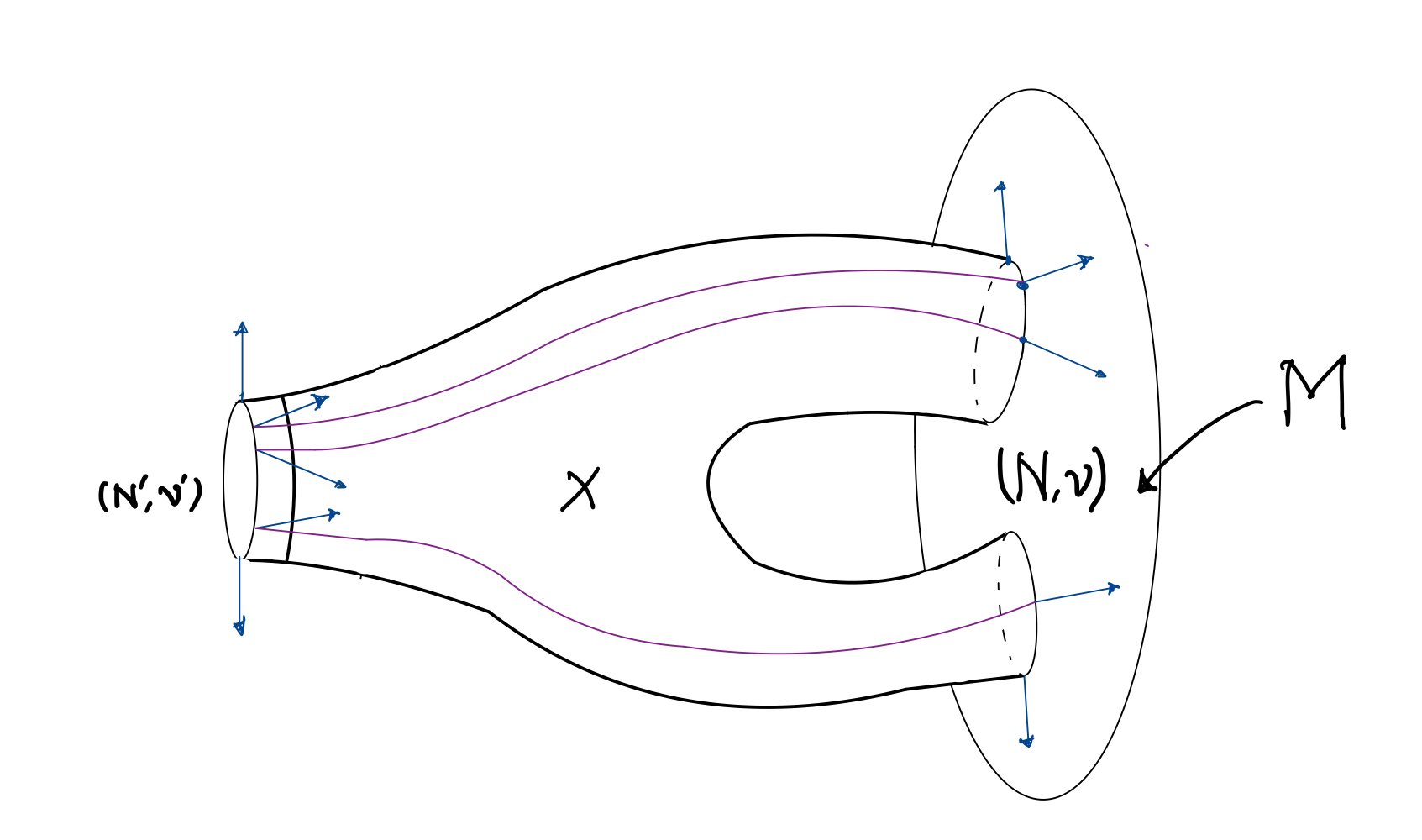}\]

\noindent Now we will introduce some terminology. Let $M$ be a manifold of dimension $n$ and $\Pi_{fr}^{p}(M)$ be the set of compact submanifolds of $M$ with codimension $p$ upto framed-cobordism. $[M,\s^p]$ is the set of all smooth maps from $M\to \s^p$ upto smooth homotopy equivalence. There is a very beautiful connection (in-fact one-one correspondence) b/w these two sets. The next few theorems will help us to get the correspondence.

\vspace*{0.3cm}

Let $f: M \to \s^p$ be a smooth map. By Sard's theorem we get a regular value $y \in \s^p$. $f^{-1}(y)$ is a co-dimension $p$ submanifold of $M$. We will construct the framing of it in the following way: If $x \in f^{-1}(y)$ then $\ker df_x = T_x (f^{-1}(y))$ and it is surjective. Choose a positively oriented basis of $T_y\s^p$ call it $w = (w^1,\cdots, w^p)$. By surejectivity of $df_x$ we can choose $\nu(x) = (\nu^1(x),\cdots,\nu^p(x)) \in (T_x(f^{-1}y)^{\perp})$ so that $\nu^i(x)$ maps to $w^i$. This gives us a framing of $f^{-1}(y)$. In other words $(f^{-1}y,\nu)$ is a framed submanifold of $M$, we denote $\mu=f^{\ast}w$. We call it Pontryagin submanifold associated to the map $f$. 

\begin{theorem} \label{thm:A}
  If $\nu$ and $\nu'$ are positively oriented basis of $T_y\s^p$ and $T_{y'}\s^p$ respectively. Then the two framed submanifold $(f^{-1}y,f^{\ast}\nu)$ and $(f^{-1}y',f^{\ast}\nu')$ are frame cobordant.
\end{theorem}

\noindent Before going to the proof of the theorem we will prove the following lemmas.

\Lem{1. If $\nu$ and $\nu'$ are positively oriented basis of $T_y\s^p$. The framed submanifold $(f^{-1}y,f^{\ast}\nu)$ and $(f^{-1}y,f^{\ast}\nu')$ are framed cobordant. }

\noindent \ts{Proof.} We know $GL_p(\R)$ has two connected components (as a topological group) and thus $GL(T_y\s^p)$ has two connected components as a topological group. Here, the conponents are determined by positve or negative determinant. Since $\mu$ and $\mu'$ are positively oriented they lie in same component. Let $gamma$ be the path between them. This gives rise to the required cobordism of $f^{-1}y \times [0,1]$.

\Lem{2. If $y$ is a regular value of $f$, and $z$ is sufficiently close to $y$, then  $f^{-1}z$ is framed cobordant to $f^{-1}y$. (Since we have seen from the previous lemma upto cobordism framed $f^{-1}y$ is unique)}

\noindent \ts{Proof.} If we consider $C$ to be the set of critical points of $f$, $f(C)$ must be closed set of $\s^p$ and thus it's compact. So there must exist $\ep$-neighborhood of $y$ contains only regular values. Choose $z$ from this $\ep$-neighborhood. Consider one parameter family (of rotations) $r_t :\s^p \to \s^p$ so that $r_0$ is identity, $r_1$ is the rotation takes $y$ to $z$ and \begin{itemize}
  \item[(i)]  $r_t$ is identity for $t \in [0,\epsilon)$ and $r_1$ for $t \in (1-\epsilon,1]$.
  \item[(ii)] each $r_t (z)$ lies on the great circle from $y$ to $z$, and hence is a  regular value of $f$. 
\end{itemize} with the help of it we can construct a homotopy $F: M\times I \to \s^p$ by $(x,t)\mapsto r_t \circ f(x)$. Not hard to see $z$ is regular value of $F$ also $F^{-1}(z) \subset  M \times I$ is a framed manifold and providing a cobordism b/w $f^{-1}(z)$ and $f^{-1}(y)$.

\Lem{3. If $f$ and $g$ are smoothly homotopic with  $y$ being the regular value for both then $f^{-1}(y)$ and $g^{-1}(y)$ are framed cobordant.}  

\noindent \ts{Proof.} Consider a homotopy $F: M \times I \to Y$ so that $F(x,t)=f(x)$ for $t \in [0,\ep)$ and $F(x,t)= g(x)$ for $t \in (1-\ep,1]$. Now we can coose a regular value of $F$, call it $z$ so that $y$ is close enough to $z$. Thus using lemma $2$ we can conclude $f^{-1}(y)$ and $g^{-1}(y)$ are framed cobordant.

\vspace*{0.2cm}

\noindent \textit{Proof of the theorem \ref{thm:A}.} Given any two point $y$ and $y'$ we can assume the frame comes from the basis $\nu$ (by lemma 1) which is positively oriented. Consider $r_t$ be the rotation so that $r_0(y)=y$ and $r_1(y') = y$. Consider the homotopy $F : M \times I \to \s^p$ given by $(x,t)\mapsto r_t(f(x))$. By the \textit{lemma 3} we can say $f^{-1}(y)$ and $f^{-1}(r_1^{-1}(y))$ are framed cobordant i.e. $f^{-1}(y)$ and $f^{-1}(y')$ are framed cobordant.

$\hfill \blacksquare$

\vspace*{0.2cm}

\noindent With the help of this lemma we can represent a Pontryagin submanifold associated to the map $f: M \to \s^p$  uniquely upto cobordism. We will represent this class of submanifolds as $\ts{Cob}_f \in \Pi^p_{fr}(M)$. Now we will prove any framed compact submanifold of $M$ is a Pontryagin manifold.
\begin{theorem} \label{thm:product}
  (\ts{Product neighborhood theorem}.) If $N$ is a compact frmasubmanifold of $M$ of codimension $p$. There is a neighborhood of $N$ in $M$ that is diffeomorphic to $N \times \R^p$. The diffeomorphism can be choosed so that $x\in N$ represent $(x,0)\in N \times \R^p$ and the normal frame $\nu (x)$ is basis of $\R^p$.
\end{theorem}

\noindent \textit{Proof.} At first, we will prove it for $M = \R^{n+p}$ (here $n$ is dimension of $N$). Consider the map $g: N \times \R^p\to M$ defined by $(x,t_1,\cdots,t_p)\mapsto x + \nu^1(x)t_1 + \cdots + \nu^p(x)t_p$. Note that $dg_{(x,0,0,..,0)}$ is invertible. Hence $g$ maps an open neighborhood of $(x,0)$ to an open set around $x \in M$ diffeomorphically. Let, $B_{\ep}$ be the open ball around $0$ of radius $\ep$. We will prove $g$ is one-one on the entire neighborhood $N \times B_{\ep}$ for some small $\ep$. If not then for every $\ep>0$ we get $(x_{\ep},t_{\ep}),(x_{\ep}',t_{\ep}')\in N \times B_{\ep}$ such that, $$g(x_{\ep},t_{\ep}) = g(x_{\ep}',t_{\ep}')$$ Since $N$ is compact $x_{\ep}\to x$ as $\ep \to 0$ and $t_{\ep}\to 0$ similarly $x'_{\ep}\to 0$ and $t'_{\ep}\to 0$ as $\ep \to 0$ but it leads to a contradiction $x_{\ep} = x'_{\ep}$. Thus there is some $\ep$ for which $f$ is one-one on $N \times B_{\ep}$, this is the neighborhood of $N$ in $M$ which is isomorphic to $N \times \R^p$ is the obvious way. Thus the statement is true for $M= \R^{n+p}$.

\vspace*{0.2cm}

For general manifold $M$ we can give it a Riemann manifold structure in the most obvious way (we will define the inner product locally and then patch the local inner products by partition of unity). So we can talk about geodesic and their length on the manifold $M$. The idea is similar consider, the map $g: N \times \R^p \to M$ given by, $(x,t)\mapsto$ the end point of the geodesic from $x$ on the direction $(t_1\nu^1(x)+ \cdots + t_p \nu^p(x))/\norm{t_1\nu^1(x)+ \cdots + t_p \nu^p(x)}$ of length $\norm{t_1\nu^1(x)+ \cdots + t_p \nu^p(x)}$ (which is unique). The rest of the part is exactly same as the above. $\hfill \blacksquare$

\begin{theorem}\label{thm:B}
  Any compact framed submanifold $N \subset M$ is a Pontryagin submanifold.
\end{theorem}

\noindent\textit{Proof.} By the previous theorem we know there is a open subset $V$ of $N$ with a diffeomorphism $\phi: V \to N \times R$ such that $\phi(N) = N \times \qty{0}$. Now consider the projection $\pi : N \times \R^p \to \R^p$. Note that $0$ is a regular value of $\pi \circ \phi$. Now we know there is an obvious diffeomorphism of $\R^p \to \s^p \setminus \qty{N}$. So consider $r$ be the map given by the composition of following maps: \[\begin{tikzcd}
	V & {N \times \R^p} & {\R^p} & {\s^p -\qty{N}} & {\s^p}
	\arrow["\phi"', from=1-1, to=1-2]
	\arrow["\simeq", from=1-1, to=1-2]
	\arrow["\pi", from=1-2, to=1-3]
	\arrow["\simeq", from=1-3, to=1-4]
	\arrow["t"', from=1-3, to=1-4]
	\arrow[hook, from=1-4, to=1-5]
\end{tikzcd}\] Here $t : \R^p \to \s^p \setminus \qty{N}$ is the diffeomorphism that sends $0$ to $S$(south pole) and $\infty$ to $N$(north pole). We can extend the map $r$ to $f:M \to \s^n$ by mapping $$f(x)=\begin{cases}
  r(x) \text{ if } x \in V \\
  \qty{N} \text{ if } x \in V^c
\end{cases}$$
Note that $f$ is a smooth function and $S \in \s^p$ is the regular value of $f$. Now note that, 
$$f^{-1}(S)=r^{-1}(S)=\phi^{-1}\circ \pi^{-1} \circ t^{-1}(S)= \phi^{-1}\circ \pi^{-1}(0) = \phi^{-1}(N \times \qty{0}) = N$$ so we are done. $\hfill \blacksquare$

\begin{theorem}\label{thm:C}
  Two maps $f, g : M \to \s^p$ are smoothly homotopic if and only if the Pontryagin manifold associated to $f$ and $g$ are framed cobordant.
\end{theorem}

\noindent \textit{Proof.} One direction is clear by lemma 3. For other direction let, $f^{-1}y$ and $g^{-1}y$ are framed cobordant with given framed cobordism $X \in M \times [0,1]$. By the previous theorem we can represent $X$ as a Pontryagin submanifold associated to $F$ by a map $F : M \times [0,1] \to \s^p$. Note that $F^{-1}_0(y) = f^{-1}y$ and $F_1^{-1}(y) = g^{-1}y$. By the following lemma we can say $F_0 \sim_{htop} f$ and $F_1 \sim_{htop} g$ so $f\sim_{htop} g$. $\hfill \blacksquare$

\Lem{4. If $f^{-1}y$ and $g^{-1}y$ are framed cobordant, $f$ and $g$ are homotopic.}

\noindent \ts{Proof.} It will be convenient to set $N=f^{-1}(y)$. The hypothesis that $f^{*} \nu=g^{*} \nu$ means that $d f_{x}=d g_{x}$ for all $x \in N$. First suppose that $f$ actually coincides with $g$ throughout an entire neighborhood $V$ of $N$. Let $h: \s^{p}-y \rightarrow \R^{p}$ be stereographic projection. Then the homotopy $$
\begin{aligned}
& F(x, t)=f(x) \quad \text { for } x \in V \\
& F(x, t)=h^{-1}[t \cdot h(f(x))+(1-t) \cdot h(g(x))] \text { for } \quad x \in M-N
\end{aligned}
$$proves that $f$ is smoothly homotopic to $g$. Thus is suffices to deform $f$ so that it coincides with $g$ in some small neighborhood of $N$, being careful not to map any new points into $y$ during the deformation. Choose a product representation $$
N \times R^{p} \rightarrow V \subset M
$$for a neighborhood $V$ of $N$, where $V$ is small enough so that $f(V)$ and $g(V)$ do not contain the antipode $\bar{y}$ of $y$. Identifying $V$ with $N \times \R^{p}$ and identifying $S^{p}-\bar{y}$ with $R^{p}$, we obtain corresponding mappings $$
F, G: N \times \R^{p} \rightarrow \R^{p}
$$ with $$
F^{-1}(0)=G^{-1}(0)=N \times 0
$$ and with $$
d F_{(x, 0)}=d G_{(x, 0)}=\left(\text { projection to } \R^{p}\right)
$$for all $x \varepsilon N$. We will first find a constant $c$ so that $$
F(x, u) \cdot u>0, \quad G(x, u) \cdot u>0
$$ for $x \varepsilon N$ and $0<\|u\|<c$. That is, the points $F(x, u)$ and $G(x, u)$ belong to the same open half-space in $R^{p}$. So the homotopy $$
(1-t) F(x, u)+t G(x, u)
$$ between $F$ and $G$ will not map any new points into 0 , at least for $\|u\|<c$. By Taylor's theorem $$
\|F(x, u)-u\| \leq c_{1}\|u\|^{2}, \text { for } \quad\|u\| \leq 1
$$ Hence, $$
|(F(x, u)-u) \cdot u| \leq c_{1} \mid\|u\|^{3}
$$ and $$
F(x, u) \cdot u \geq\|u\|^{2}-c_{1}\|u\|^{3}>0
$$ for $0<\|u\|<\operatorname{Min}\left(c_{1}^{-1}, 1\right)$, with a similar inequality for $G$. To avoid moving distant points we select a smooth map $\lambda: R^{p} \rightarrow R$ with $$
\begin{array}{lll}
\lambda(u)=1 & \text { for } & \|u\| \leq c / 2 \\
\lambda(u)=0 & \text { for } & \|u\| \geq c
\end{array}
$$ Now the homotopy $$
F_{t}(x, u)=[1-\lambda(u) t] F(x, u)+\lambda(u) t G(x, u)
$$ deforms $F=F_{0}$ into a mapping $F_{1}$ that (1) coincides with $G$ in the region $\|u\|<c / 2$, (2) coincides with $F$ for $\|u\| \geq c$, and (3) has no new zeros. Making a corresponding deformation of the original mapping $f$, this clearly completes the proof of Lemma 4. $\hfill \blacksquare$

\vspace*{0.2cm}

\noindent With the help of theorem \ref{thm:A},\ref{thm:B},\ref{thm:C} we can conclude $\Pi^p_{fr}(M) = [M,\s^p]$ as a set (infact the later can be given a group structure discussed below). This is called \textit{cohomotopy group}.

\vspace*{0.2cm}

\noindent Let $m$ be the dimension of $M$. We can give $\Pi_{fr}^p(M)$ a group structure for certain $p$'s. If $N$ and $N'$ are two submanifold of $M$ of codimension $p$, then their transversal intersection is also a submanifold of codimension $2p$. We want the transversal intersection to be empty (so that we can give disjoint union a group operation). In other-words $\sim N + \dim N' < \dim M$ thus we have $m-p+m-p \leq m-1$ and thus $p \geq \frac{1}{2}m+1$. Now the operation $\sqcup : \Pi_{fr}^p(M) \times \Pi_{fr}^p(M)\to \Pi_{fr}^p(M)$ gives a product structure on $\Pi^p_{fr}(M)$ in-fact it is an Abelian group. The identity element of this group is the class $[\emptyset]$ consisting of all closed submanifolds which are boundaries of some manifold. For any manifold $(M,\nu)$ with it's framing we can consider $(M,-\nu)$ (the opposite framing), if we denote the former by $M$ and the later by $[-M]$. One can check $[M]+[-M]=[\emptyset]$. Now we can define a product 

$$\otimes : \Pi^p_{fr}(M) \times \Pi^q_{fr}(M)\to \Pi^{p+q}(M)$$

\noindent Which is given by transversal intersection. If $N,N'$ are submanifold of codimension $p$ and $q$ respectively we can perturb $N$ so that $N$ and $N'$ intersect transversally. For transversal intersection $N \cap N'$ is a submanifold of codimension $p+q$ (details can be found in \textcolor{gray}{Guillemin and Pollack, ch1}). Thus we can define a graded ring $$\Pi^{\ast}_{fr} (M) = \bigoplus_{p \geq \frac{1}{2}m+1} \Pi^p_{fr}(M)$$

\section{Hopf's Theorem and $\pi_n(\s^n)$ or $\pi^{S}_0$}

If $M$ is oriented connected and boundaryless manifold of dimension $m$.  A framed submanifold of codimension $m$ is given by $f^{-1}(p),f^{\ast}\nu$ for some smooth map $f :M \to \s^{m}$. Now this $f^{-1}$ is nothing but finite set of points with the subspace topology is the discrete topology. The $f^{\ast}\nu$ will have some orientation for $f^{-1}p = \qty{q_1,\cdots, q_r}$. We associate $+1$ is $f^{\ast}(\nu)(q_i)$ have same orientation as $\nu$ otherwise we will associate $-1$. We denote this by $\text{sgn}(q_i)$. It's not difficult to see that the framed cobordism class of $0$-dimensional submanifold of $M$ is uniquely determined by $\sum \text{sgn}(q_i)$. Now the sum $$\sum_{i} \text{sgn}(q_i) = \deg (f)$$ so we can conclude the following theorem

\begin{theorem}
  \ts{Hopf's theorem} If $M$ is a connected, oriented and boundaryless manifold of dimension $n$ two maps $M \to \s^n$ are homotopic iff their degree is same.
\end{theorem}

\noindent Given any integer $n \in \Z$ we can construct a map $f: \s^m \to \s^m$ which have degree $n$. Using Hopf's theorem we can say $\pi_n(\s^n)=[\s^n,\s^n] =\Z$.

\vspace*{0.2cm}

\noindent \ts{Remark.} The above theorem is not true for un-oriented manifolds. But if we look at degree mod $2$. Then the above theorem is still true. We sometime use the fact $\Pi_{fr}^k(\R^{n+k}) = \pi_{n+k}(\s^k)$ to compute the higher homotopy groups of sphere. Now we define $\Pi_n^{fr}(\R^{n+k}):= \Pi_{fr}^k(\R^{n+k})$, in other-words it's the framed cobordism class of $n$ dimensional smooth submanifolds of $\R^{n+k}$. The framing of a manifold $M \subset \R^N$ exist if the normal bundle of $M$ is trivial, so We can treat $M$ as a submanifold of $\R^{N+1}$, here the normal bundle is also trivial and isomorphic ro $N_{\R^N}(M)\oplus \ep$. So if $f$ is framing of $M \subset \R^N$ there is a natural framing of $M \subset \R^{N+1}$ by the natural inclusion of vectors in $f(x)$ in side $\R^{N-\dim M +1}$ and the last vector being $(0,0,\cdots,1)$; call it $f \oplus \ep$. So, there is a natural inclusion $\iota: \Pi_n^{fr}(\R^{n+k}) \hookrightarrow \Pi_n^{fr}(\R^{n+k+1})$ by $[(M,f)]\mapsto [(M,f \oplus \epsilon)]$ and we can show the following diagram commutes: \[\begin{tikzcd}
	{\Pi_n^{fr}(\R^{n+k}) } & \Pi_n^{fr}(\R^{n+k+1})  \\
	{\pi_{n+k}(\s^k)} & \pi_{n+k+1}(\s^{k+1})
	\arrow["\iota", from=1-1, to=1-2]
	\arrow["{\text{Thom map}}"', from=1-1, to=2-1]
	\arrow["\text{Thom map}", from=1-2, to=2-2]
	\arrow["-\wedge 1"', from=2-1, to=2-2]
\end{tikzcd}\] \label{comm:1} We can take colimit and it will give us : $$\underset{k}{\colim} \, \Pi_n^{fr}(\R^{n+k}) \simeq \pi_k^S$$ If we define $\Omega_n^{fr}$ to be the set of all smooth $n$-dim manifold quotiented by the equivalance relation induced from framed cobordism. Note that any $n$-smooth manifold $M$ must lie inside some $\R^{n+k'}$ so from the embedding we get a framing of the manifold $M$, and so, $M \in \Pi_n^{fr}(\R^{n+k'})$ for the choice of $k'$. And thus we can say $$\Omega_n^{fr} = \underset{k}{\colim} \, \Pi_n^{fr}(\R^{n+k}) \simeq \pi_n^S$$  


\noindent With the help of the theories developed above, we will compute sable homotopy group of spheres for a few indices. 

\section{The first stem: $\pi^S_1$}

We begin with the Hopf fibration $\eta : \s^3 \to \s^2$ with the homotopy fiber $\s^1$. From the \textit{Puppe sequence}, we deduce that $\pi_n(\s^2) \simeq \pi_n(\s^3)$ for $n \geq 3$, and hence $\pi_3(\s^2) \simeq \pi_3(\s^3)$. The generator of this group is the homotopy class of $\eta$. On the other hand, from [\ref{thm:freu}], we know that $\pi_{n+1}(\s^n)$ stabilizes for $n \geq 3$. To calculate $\pi_1^S$, it suffices to compute $\pi_4(\s^3)$. The Freudenthal suspension theorem tells us that the map $\Sigma : \pi_3(\s^2) \to \pi_4(\s^3)$ is surjective. Thus, $\pi_4(\s^3)$ must be generated by $[\Sigma \eta]$. Our goal in this section is to analyze $\Sigma \eta$ and its relation to determining $\pi_1^S$. 

\vspace{0.3cm}

\noindent Here, we use the identification of $\Pi_n^{fr}(\R^{n+k})$ with $\pi_{n+k}(\s^k)$. There exists a framed $1$-manifold in $\s^3$ (or $\R^3$) corresponding to the Hopf map $\eta$. From the Thom-Pontryagin construction, we know that this manifold can be described as the inverse image of a regular value. Thus, it must be $\s^1$, as the fiber of $\eta$ is $\s^1$ for any chosen point in $\s^2$. From the following commutative diagram, we can choose a regular value of $\Sigma \eta$ such that $(\Sigma \eta)^{-1}(p) = (\eta \times \text{Id})^{-1}(p)$. This inverse image should be homeomorphic to $\s^1$: 

\[\begin{tikzcd}[cramped]
	{\s^3\times \s^1} & {\s^2\times \s^1} \\ 
	{\s^4} & {\s^3} 
	\arrow["{\eta \times \text{Id}}", from=1-1, to=1-2] 
	\arrow[from=1-1, to=2-1] 
	\arrow[from=1-2, to=2-2] 
	\arrow["{\Sigma \eta}"', dashed, from=2-1, to=2-2] 
\end{tikzcd}\]

\noindent Therefore, the Pontryagin manifold in $\Pi_1^{fr}(\R^4)$ corresponding to $[\Sigma \eta]$ is a circle embedded in $\R^4$. To describe the homotopy class $\Sigma \eta$, we need to consider the framing of this Pontryagin manifold. From the commutative diagram [\ref{comm:1}], it suffices to determine the framing of the Pontryagin manifold corresponding to $\eta$. 

\vspace*{0.2cm}

\noindent The following diagram illustrates two types of framings of a circle. The first type represents $0$ in $\Pi_1^{fr}(\s^3)$ since it can be extended to a framing of a disk. However, the second type of framing corresponds to a non-zero element in $\Pi_1^{fr}(\s^3)$ (recall that $\Pi_1^{fr}(\s^3) \simeq \Z$). 

\begin{figure}[h]
    \centering
    \includegraphics[width=0.6\textwidth]{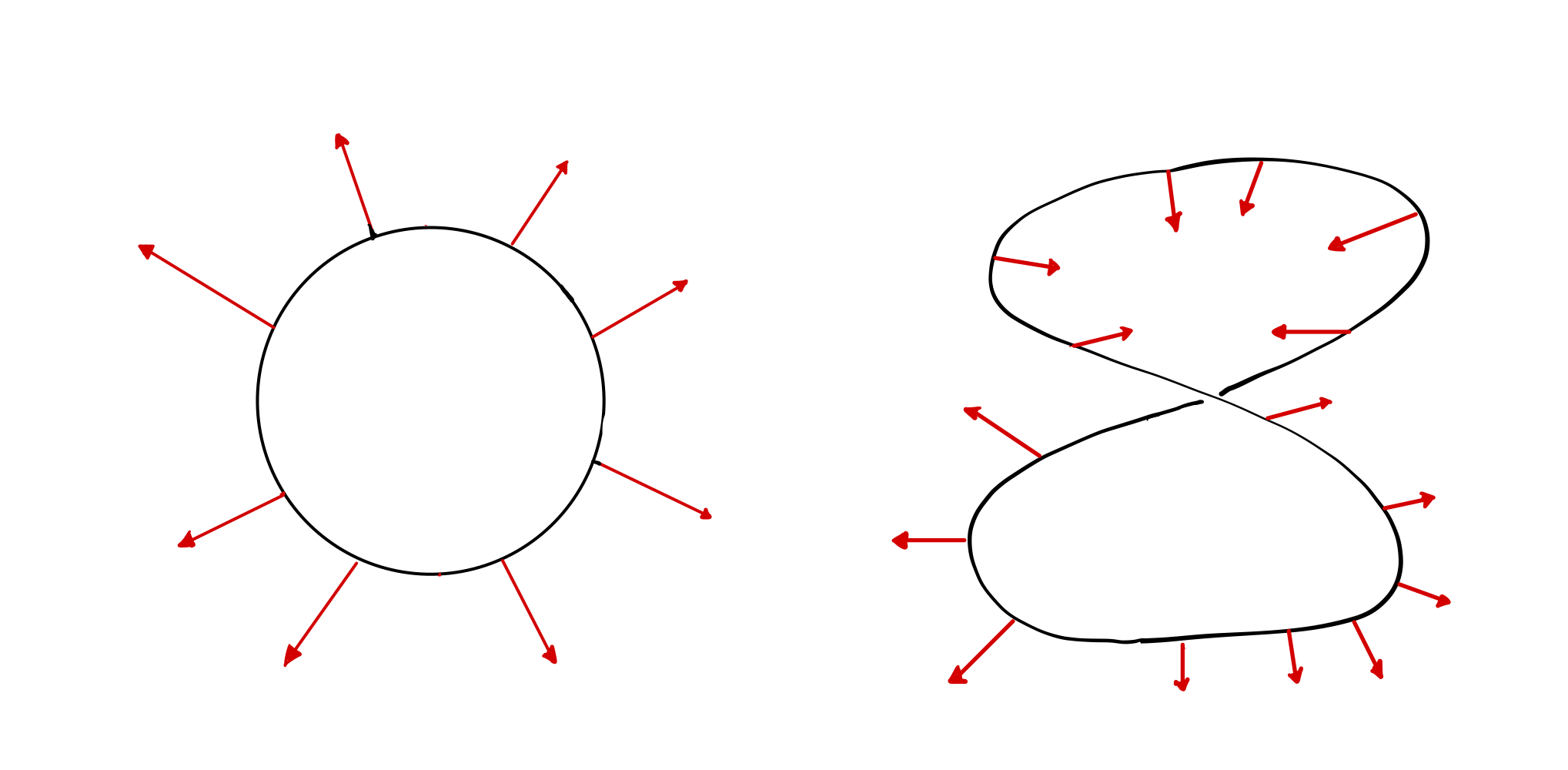}
    \caption{\textcolor{purple}{Possible framings of $\s^1\hookrightarrow \R^3$}}
    \label{fig:mesh1}
\end{figure}

\noindent We will provide a proof, inspired by \cite{Pont1}, showing that the Pontryagin manifold corresponding to $\Sigma \eta$ is associated with the second type of framing. 

\vspace*{0.2cm}

\noindent Let, $(M,\nu)$ be a framed $1$-manifold(compact) in $\R^4$. Thus by classification of $1$ manifolds $M$ is disjoint union of $\s^1$ (upto framed cobordism). Let, $\nu(x)$ is positively oriented (orientation cominng from the tanget bundle) $\forall x \in M$. Note that $\nu(x) = (\nu^1(x),\cdots,\nu^3(x))$ is a positively oriented basis of $N_x(M \subset \R^4)$. Define $\nu'(x)$ be the element of normal-space $N_x(M \subset \R^4)$, we get after Gram-Schmidt orthonormalization of $\nu(x)$. The deformation retract of $GL_n(\R)$ to $O(n)$, will help us to give us a framed cobordism between $(N,\nu)$ and $(N, \nu')$. Without loss of generality we may assume $(M,\nu)$ is a $1$-manifold in $\R^4$ with $\nu(x)$ is positively oriented orthonormal basis of the normal space. For every $x\in M$ there is a unique vector $\tau(x)$ in $T_xM$ so that $(\tau(x),\nu(x))$ is an element of $SO(4)$ (with respect to the standard basis). We can define a map $$h^{M,\nu}: M \to SO(4)$$ given by $x \mapsto (\tau(x), \nu(x))$; clearly this is a continuous map. Let, $[M]$ be the fundamental class of $M$ and then $h_{\ast}^{M,\nu}([M]) \in H_1(SO(4);\Z)$ is an element of $\Z/2\Z$ $^{\textcolor{red}{\ast \ast}}$, We define residue class of $(M,\nu)$ by \newcommand{\res}{\textbf{Res}} $$\res(M,\nu) = h_{\ast}^{M,\nu}([M]) + \text{ no of components in } M \, \pmod{2}.$$  The above definition of residue is given for the standard orientation on $\R^4$; this is independent of orientation on $\R^4$, if $(\tau(x), \nu(x))$ is not positively oriented we can take $(-\tau(x), \nu(x))$ to be positively oriented basis (i.e determinant w.r.t basis is $+1$). But $(\tau(x),\nu(x))\mapsto (-\tau(x),\nu(x))$ is homeomorphism so the image of fundamental class under $h_{\ast}^{M,\nu}$ doesn't change upto sign. So residue is independent of orientation on $\R^4$.

\begin{lemma}{}{} \label{lem-res:1}
	\hspace{0.2cm} If two framed manifold $(M_0,\nu_0)$ and $(M_1,\nu_1)$ are framed cobordent then $$\res(M_0,\nu_0)=\res(M_1,\nu_1)$$
\end{lemma}

\noindent In order to prove the above proposition, we recall some results from Morse theory and low dimensional topology. Let, $f: M \to \R$ be a smooth function; a point $p$ is said to be \textit{critical} if $\nabla f(p) =0$ also it is said to be \textit{degenerate} if $\nabla f (p)=0 = \det H_f (p)$ (here $H_f(p)$ is Hessian). A function is said to be \textit{Morse function} if it do not have any degenerate critical function. The following are important results from Morse theory \cite{Morse1} we will be using here:

\begin{itemize}
	\item[-] If $f:M \to \R$ is a smooth function from a compact manifold, it can be approximated arbitrarily by a Morse function. 
	\item[-] Any smooth function around a non-degenerate critical point can be written as $f = -(x_1^2+\cdots +x_k^2) + (x_{k+1}^2+\cdots + x_n^2)$ with respect to a coordinate chart around $p$ in $M$.  Here $n= \dim M$ and $k$ is the number of negative eigenvalues of Hessian.
\end{itemize}

\noindent In low dimension topology we heavily use \textit{Handle body} decomposition. A $n$-dimensional $k$-handle is the manifold $D^k \times D^{n-k}$. By attaching a $k$-handle we mean attaching $D^k \times D^{n-k}$ along $\p D^k \times D^{n-k}$.

\begin{theorem}
	(Handlebody decomposition from Morse function) Let, $f:M \to \R$ be a morse fnction and $[a,b]\subset \R$ be an interval where we have only one critical point. Then, The manifold $f^{-1}(-\infty,b]$ can be achived by attaching $k$-handle to $f^{-1}(-\infty,a]$. Where $k$ is the index at the critical point. 
\end{theorem}

\noindent The proof can be found in \cite{Morse1}.

\begin{figure}[h]
    \centering
    \includegraphics[width=0.9\textwidth]{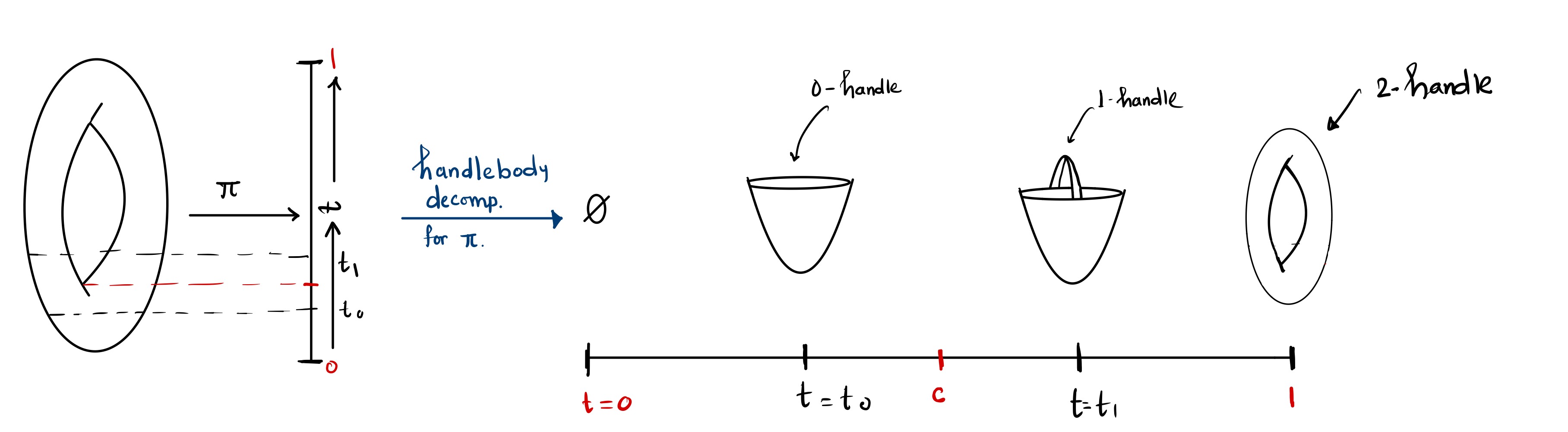}
    \caption{\textcolor{purple}{Example of Handlebody decomposition of a morse function}}.
    \label{fig:mesh2}
\end{figure}

\noindent \textbf{Proof of the proposition.}[\ref{lem-res:1}] Let, $M$ be the framed cobordism of the manifolds $M_0,M_1$, $M \subset \R^4 \times [0,1]$. There is a smooth projection map $\pi: M \to [0,1]$ by restricting the natural projection $\R^4 \times [0,1]$ to $[0,1]$. By the Morse approximation let, $\pi$ be morse function. Let, there is no critical value of $\pi$ on $[0,\ep]$ then $\pi^{-1}(0)=M_0$ and $\pi^{-1}(\ep)=M_{\ep}$ are diffeomorphic, so $\res(M_0,\nu)= \res(M_{\ep},\nu_{\ep})$ here, $\nu_{\ep}$ is the framing induced from the framing of $M$. Let, $M_t = \pi^{-1}(t)$ and $\nu_t$ is the framing induced from $M$, for regular value $t$. The residue value can chanage only if we pass through a critical point. 

\vspace{0.2cm} \newcommand{\de}{\delta}

\noindent Let $c$ be a critical value of $\pi$, we know for morse function critical values are isolated, so we get a neighborhood $[c-\de,c+\de]$ so that it has only one critical value $c$. Let, $M^{-}= \pi^{-1}[0,c-\de]$ and $M^{+}= \pi^{-1}[0,c+\de]$, If we show the residue are same for $M_{c-\de}$ and $M_{c+\de}$ we are done. 

\vspace{0.2cm}

\noindent By the Handlebody decomposition theorem we can say $M^{+}$ can be achived by attaching a $k$-handle to $M^{-}$. If we aattach a $0$-handle, in that case We are adding a adittional component $C$ to $M_{c-\de}$ to get $M_{c+\de}$. The component $C$ encloses a framed disk, we can treat it like a trivially framed disk in $\R^2$ thus $\res(C,\nu_C)=0$. So $\res(M_{c-\de},\nu_{c-\de})= \res(M_{c+\de},\nu_{c+\de})$.

\vspace{0.2cm}

\noindent If we attach a $2$-handle to $M^{-}$ then also, $M_{c+\de}$ can be achived by adding a component $C$ that encloses a framed disk to $M_{c-\de}$ or by attaching $\s^1 \subset D^2$ to $M_{c-\de}$. The former case is similar to the $0$-handle attaching. The later case is also similar as $M_{c+\de}$ is a framed circle that encloses a framed disk.

\vspace{0.2cm}

\noindent We are left with the case when we attach $1$-handle to $M^{-}$ to get $M^{+}$. Let, $\pi(m)=c$, there is a co-ordinate chart around $m$ where $\pi$ looks like $x_1^2-x_2^2$ and the co-ordinate of $m$ is $(0,0)$. This is a cross, the component of $M_{c}$ containing this cross must be a `figure eight' space; call it $L$. For small $\de$ and $c+\de>t >c$; the part of $M_t$ near $L$ is made of one circle $C_0$ and for $c-\de<t< c$; the part of $M_t$ near $L$ is made of two circle $C_1,C_2$. Let the induced framing of $C_0,C_1,C_2$ be $\nu_0,\cdots,\nu_2$ respectively. Our aim is to show $$h_{\ast}^{C_0,\nu_0}([C_0])+1 = h_{\ast}^{C_1,\nu_1}([C_1])+ h_{\ast}^{C_2,\nu_2}([C_2]) \pmod{2}$$ 

\noindent -------------------------------------------------------------------------------------------------------------------------------------

\noindent $\textcolor{red}{\ast \ast}$ This is because the fundamental group of topological groups are abelian and since $SO(4)$ is path-conneceted, by Hurewicz Theorem the homology group should be isomorphic to $\pi_1(SO(4))$, the CW -decomposition of $SO(4)$ will give us that it's $2$-skeleta is $SO(3)\simeq \R P^3$, so $H_1(SO(4);\Z) \simeq \Z/2\Z$. 

\noindent -------------------------------------------------------------------------------------------------------------------------------------

\noindent $\textcolor{red}{\ast \ast}$ The $\res$ can be defined for any embedding $M \hookrightarrow \R^n$.

\[
    \includegraphics[width=0.85\textwidth]{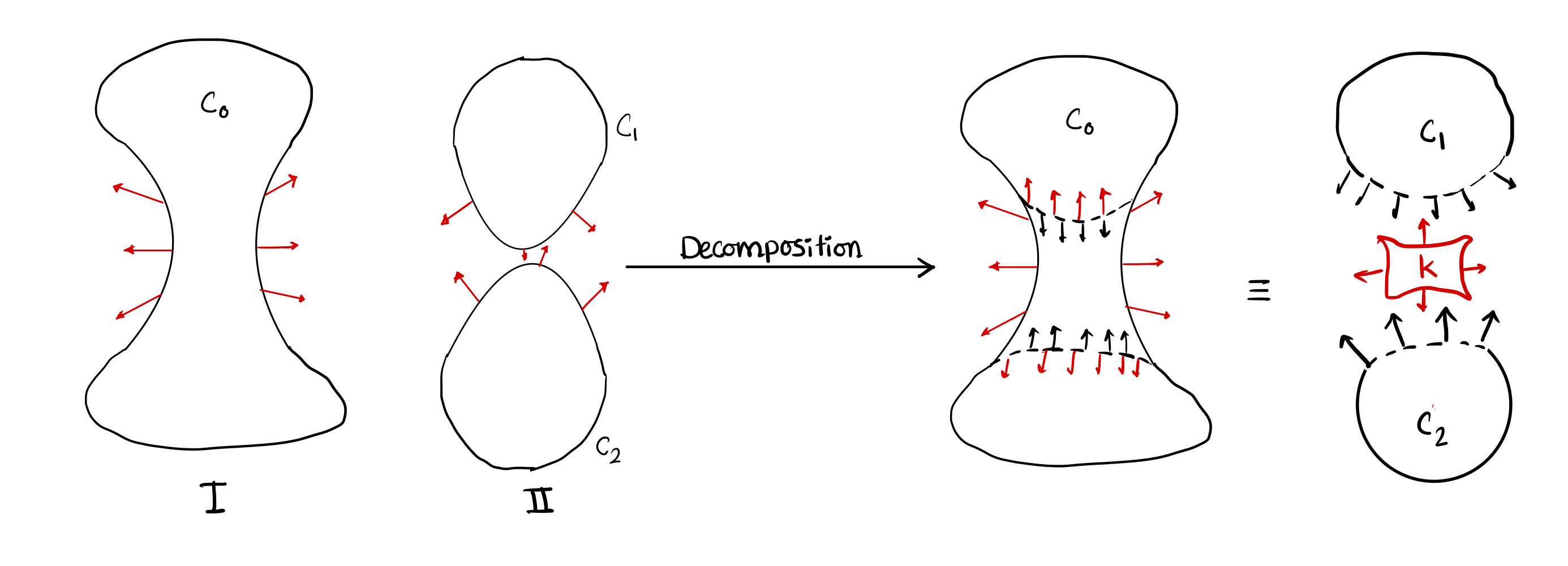}
\]    \label{fig9}

\noindent  The above picture shows us the framing on $C_0, C_1,C_2$ (locally). For the framing we have used the local description of $\pi$'s. For (I), it's locally given by $x_1^2-x_2^2=c'>0$, which is a hyperbola and the framing comes from $\nabla \pi$; similarly we got the framing on $C_1,C_2$. So, we can use the decomposition of $C_0$ in to $C_1,C_2$ and $K$ where the framing on $K$ is the standered framing of circle in $\R^2$, thus $$h_{\ast}^{C_0,\nu_0}([C_0])+h_{\ast}^{K,\nu_K}([K]) = h_{\ast}^{C_1,\nu_1}([C_1])+ h_{\ast}^{C_2,\nu_2}([C_2]) \pmod{2}$$ which completes the proof. $\hfill \blacksquare$

\vspace*{0.4cm}

\noindent Let, we have a continuous map $g: M \to SO(3)$ (here we are considering $M$ to be a framed manifold with the framing $\nu$). This gives another framing on $M$; With respect to sstanderd co-ordinate if we can write; $g(x) = (g_{ij}(x))$ then the new framing on $M$ can be given by $$x \mapsto (\sum_j g_{ij}(x)\nu^j(x))_{i=1}^3$$ we denote it by $g(\nu)$ and call this a twist of the framing $\nu$. Let, $\s^1 \subset \R^4$ be the circle with the framing comes from $\R^2$ (call the canonical framing $xi$); then for any continuous map $g: \s^1 \to SO(3)$ we have a framing of $\s^1$. If two such loops are homotopic, then we can use the homotopy to construct a framed cobordism between corresponding framed circles. Thus we can define a map a follows: $$J_3: [\s^1, SO(3)]\to \Pi_{1}^{fr}(\s^4)$$ Indeed we can carry out the same work for $\s^1 \subset \R^{n+1}$ and define $J_n:[\s^1,SO(n)]\to \Pi_1^{fr}(\s^{n+1})$. Note that $J_2$ is isomorphism. From the commutativity of [\ref{comm:1}] we can say $\iota : \Pi^{fr}_1(\s^3) \to \Pi^{fr}_1(\s^{n+1})$ is surjective as $\Sigma^{n-2}$ is. There is a inclussion of $i: SO(2) \hookrightarrow SO(n)$ given by $A \mapsto \text{diag}[A,I_{n-2}]$. Thus the following diagram commutes.   \[\begin{tikzcd}[cramped]
	{\pi_1(SO(2))} & {\Pi^{fr}_1(\s^{3})} \\
	{\pi_1(SO(n))} & {\Pi^{fr}_1(\s^{n+1})}
	\arrow["{J_2}", from=1-1, to=1-2]
	\arrow["{i_{\ast}}"', from=1-1, to=2-1]
	\arrow["\iota", from=1-2, to=2-2]
	\arrow["{J_n}"', from=2-1, to=2-2]
\end{tikzcd}\] Here $J_2$ is isomorphism thus $J_n$ is surjective. We know $\pi_1(SO(n)) \simeq \Z/2\Z$, thus it has a generator. Let, the class of $f: \s^1 \to SO(n)$ be the generator, then $(\s^1,f(\xi))$ is a framed $1$-manifold of $\s^{n+1}$. Now we claim that $\res(\s^1,f(\xi))=1$ and so it's a nontivial element in the cobordism class (when $n=3$), so $J_n$ is an isomorphism. Note that, $$h^{\s^1,f(\xi)}(x)= \begin{pmatrix}
	1 & 0 \\
	0 & f(x) 
\end{pmatrix}h^{\s^1,\xi}(x)$$ Now note that the map $\text{diag}[1,f] : \s^1 \to SO(n+1)$ gives us the map same with $f_{\ast}: H_1(\s^1) \to H_1(SO(3))$ in homology. If we have two maps $\alpha, \beta : \s^1 \to SO (3)$ then $(\alpha \beta)_{\ast}= \alpha_{\ast} + \beta_{\ast}$ (\cite{hatat};{chapter 1}). Thus $$\res(\s^1,f(\xi))= h^{\s^1,f(\xi)}_{\ast}([\s^1])+1 = f_{\ast}([\s^1]) + h^{\s^1,\xi}_{\ast}([\s^1])+1= f_{\ast}([\s^1]=1)$$ The last equality is because of Hurewicz isomorphism. Thus $\Pi_1^{fr}(\s^4) \simeq \Z/2\Z$. From the correspondence with homotopy group we can conclude $[\Sigma \eta]$ is non-zero element in $\pi_3(\s^4)$ and it has order $2$. Thus $$\pi_3(\s^4)\simeq \Z/2\Z \simeq \pi_1^S$$

\section{The second stem : $\pi_2^S$}

In order to compute the stable stem \(\pi_2^S\), we need to determine \(\pi_3(\s^5)\). From the Hopf fibration, we already have \(\pi_4(\s^2) \simeq \pi_4(\s^3)\), with the generator of \(\pi_4(\s^2)\) given by \([\eta \circ \Sigma\eta]\). By the Freudenthal suspension theorem, we know that the suspension map \(\Sigma : \pi_4(\s^2) \to \pi_5(\s^3)\) is surjective. Consequently, the generator of \(\pi_5(\s^3)\) is \([\Sigma (\eta \circ \Sigma \eta)]\). We will develop tools to show that \(\alpha := [\Sigma (\eta \circ \Sigma \eta)]\) represents a nontrivial homotopy class. This will establish that the above map is an isomorphism, leading to the conclusion that \(\pi_2^S \simeq \mathbb{Z}/2\mathbb{Z}\).

\vspace{0.2cm}  

\noindent There are several approaches to proving the nontriviality of \(\alpha \in \pi_5(\s^3)\). Modern techniques use Steenrod operations, providing an algebraic and relatively straightforward proof. However, the aim of this article is to present a geometric proof. Thus, we utilize the most direct tool available: the correspondence between framed cobordant 2-manifolds in \(\mathbb{R}^5\) or \(\mathbb{S}^5\) and elements of \(\pi_5(\s^3)\). The Pontryagin manifold corresponding to the map \(\Sigma(\eta \circ \Sigma \eta)\) is a torus with an obvious framing. Our task is to show that this framed torus is not null-cobordant in \(\mathbb{R}^5\).

\vspace{0.2cm}

\noindent As in the previous section, we consider orthonormal framings of manifolds and aim to develop a homology invariant to carry out similar computations. We have already shown that ordinary cohomology theories are representable by the Eilenberg–MacLane spectrum. Thus, for an orientable compact 2-manifold \( M \), every cycle in \( H_1(M;\mathbb{Z}) \) corresponds to a cocycle in \( H^1(M;\mathbb{Z}) = [M,\mathbb{S}^1] \), which establishes a one-to-one correspondence with Pontryagin’s 1-dimensional manifolds in \( M \).  

\vspace{0.2cm}

\noindent We can associate each homology class of a cycle with a framed 1-submanifold of \( M \), yielding a natural map  
\[
\hat{q}_{M}: H_1(M;\mathbb{Z}) \to \mathbb{Z}/2\mathbb{Z}
\]
defined by \( x \mapsto \res(x) \).  For \( x, y \in H_1(M;\mathbb{Z}) \), we define the corresponding submanifold in the same manner. The sum \( (x + y) \) represents the disjoint union of \( x \) and \( y \) with their framings. By perturbing these submanifolds, we ensure that they intersect transversally. Since \( x \) and \( y \) are compact, their intersection \( x \cap y \) is finite. We define the intersection number  
\[
I(x,y) = \text{number of intersections of } x, y \pmod{2}.
\]  
Let \( \{p_1, \dots, p_m\} \) be the set of intersection points. Consider \( K \) as the union of \( x \) and \( y \) such that it contains these intersection points. The manifold \( K \) inherits a framing from \( M \), leading to the relation  
\[
\res(K) = \res(x) + \res(y).
\]  
By applying the surgery described in Figure [\ref{fig9}], we obtain \( x \sqcup y \). At each point \( p_i \), an additional component of \( K \) appears, yielding  
\[
\res(x+y) = \res(x) + \res(y) + I(x,y) \pmod{2}.
\]  
Clearly $\hat{q}_M$ is not a homomorphism. In the paper \cite{Pont2} this map was proven to be a homomorphism, which is not true by above discussion. Also, it is straightforward to verify that \( \hat{q}_M(x + x) = 0 \), allowing us to define a map  
\[
q_M: H_1(M;\mathbb{Z}/2\mathbb{Z}) \to \mathbb{Z}/2\mathbb{Z}.
\]  
We call such functions \textit{quadratic refinements} of the symmetric bilinear form \( I(-,-) \).

\vspace{0.2cm} \newcommand{\arf}{\mathbf{Arf}}

\noindent Let \((M,\nu)\) be a connected, compact, orientable framed manifold in \(\mathbb{R}^5\) (although the specific dimension 5 is not essential for the following definitions, it will become relevant later in our discussion). By the classification of surfaces, we can define the genus of \( M \), denoted by \( g(M) \).  

The homology group \( H_1(M;\mathbb{Z}/2\mathbb{Z}) \) forms a vector space of dimension \( 2g(M) \) over \( \mathbb{Z}/2\mathbb{Z} \). Since \( I \) is a bilinear form, there exists a basis \( \{a_1, \dots, a_{g(M)}, b_1, \dots, b_{g(M)}\} \) such that  
\[
I(a_i, b_i) = 1, \quad I(a_i, a_j) = I(b_i, b_j) = I(a_i, b_j) = 0 \text{ for } i \neq j.
\]  
As discussed earlier, we may consider \( a_i, b_i \) as framed submanifolds of \( M \). We define the *Arf invariant* as  
\[
\arf(M,\nu) = \sum_{i=1}^{g(M)} q_M(a_i) q_M(b_i) \pmod{2}.
\]  
Notably, this definition is independent of the choice of basis. Although Pontryagin did not explicitly refer to this quantity as the "Arf invariant" in his original work, it was later recognized that it also appeared in Arf’s research.  

\vspace{0.2cm}

\noindent We will now prove that the Arf invariant is a framed cobordism invariant and show that for the Pontryagin manifold corresponding to \( \alpha \), this invariant is \( 1 \). This establishes that \( \alpha \) represents a non-trivial class in \( \pi_5(\mathbb{S}^3) \), completing our proof.

\begin{theorem}\label{arf:1}

	\hspace{0.1cm} If two framed $2$-manifold $(M_0,\nu_0)$ and $(M_1,\nu_1)$ in $\R^5$ are framed cobordent then $$\arf(M_0,\nu_0)=\arf(M_1,\nu_1)$$
\end{theorem}

\noindent \textbf{Proof.} The proof will be very similar to the proof of [\ref{lem-res:1}]. Let, $M\subset \R^5 \times [0,1]$ be the framed cobordism between $M_0,M_1$. Upto isotopy we can assume the restriction of the projection $\pi: \R^5 \times [0,1]\to [0,1]$ from $M \times [0,1]$ is a Morse function. Define $M_t = \pi^{-1}(t)$ it's clear that $\arf(M_t)$ can change only if $t$ pass through a critical value. Let, $c \in [0,1]$ be the minimum critical value of $\pi$, since critical values are isolated we may assume tere is a neighborhood $[c-\de,c+\de]$ which has only critical value. Our proof will be done if we can show $\arf(M_{c-\de})=\arf(M_{c+\de})$ (we are always keepig track of the framing by taking induced framing from $M$). Let, $M^{+}=\pi^{-1}[0,c+\de], M^{-}=\pi^{-1}[0,c-\de]$. Similar to [\ref{lem-res:1}] we have four cases here: 

\vspace{0.2cm}

\noindent If we get $M^{+}$ by attaching $0$-handle to $M^{-}$. In this case $c$ is a local minima and $M_{c+\de}$ differs from $M_{c+\de}$ by one $D^3$ which has trivial homology so $\arf$ remains unchanged.  Similarly for $3$-handle attachment $\arf$ remians unchanged.

\vspace{0.2cm}

\noindent If we get $M^{+}$ by attaching $1$-handle to $M^{-}$. We basically add a tube to $M_{c-\de}$ to get $M_{c+\de}$, the framing on these manifolds can be given by the framing on $M$. Note that the genus get increased by $1$ in this process. 

\[
    \includegraphics[width=0.65\textwidth]{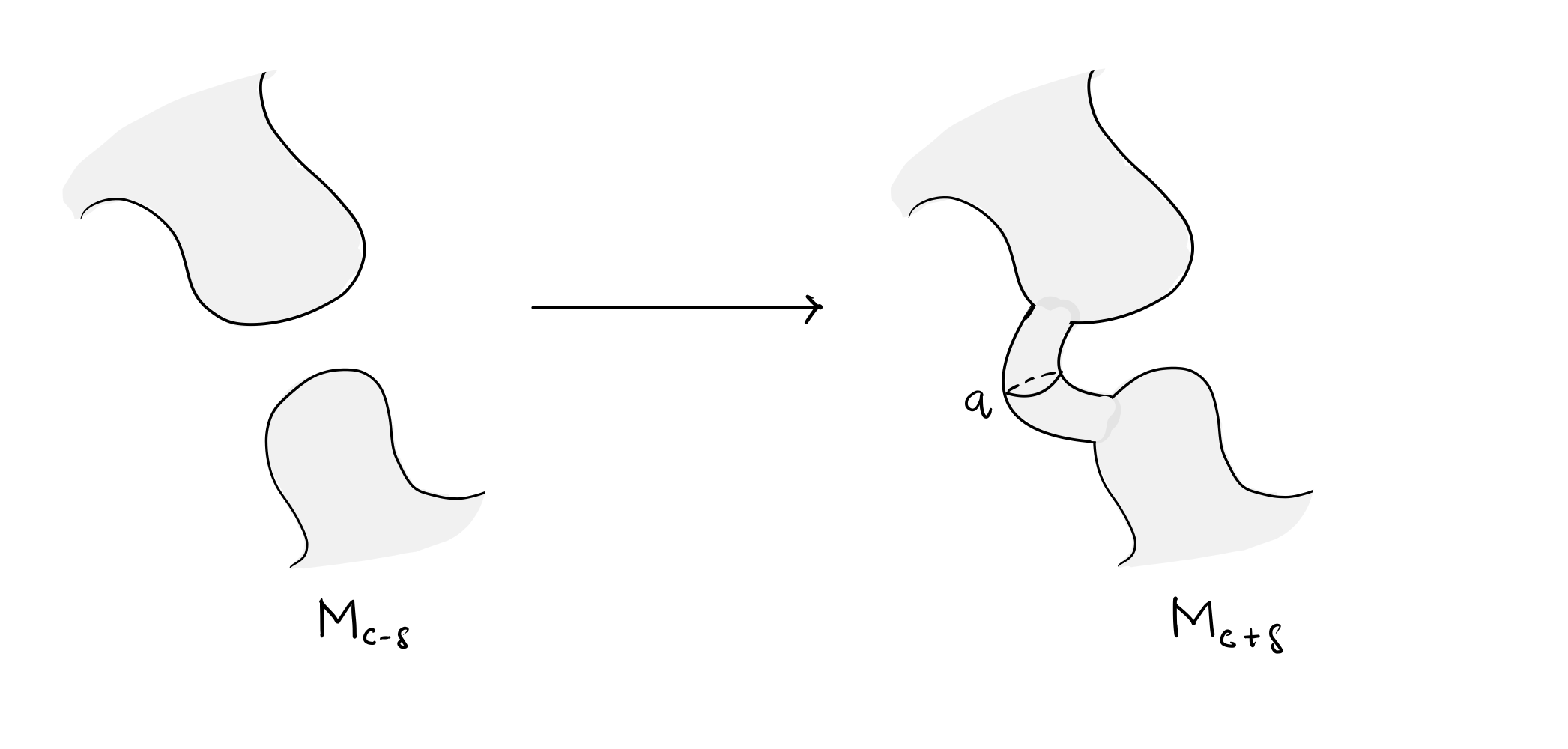}
\]    \label{fig10} Let,$g=g(M_{c-\de})$ Note that the submanifold $a$ bounda a disk in $M$ so $\res(a)=0$, we can generate a basis of $H_1(M_{c+\de};\Z/2\Z)$, $\qty{a_i,b_i}_{i=1}^{g+1}$ such that $a_{g+1}=a$ and $\qty{a_i,b_i}_{i=1}^{g}$ are in $M_{c-\de}$. So $$\arf(M_{c+\de})=\sum_{i=1}^{g}q_{M_{c-\de}}(a_i)q_{M_{c-\de}}(b_i)+0 = \arf(M_{c-\de})$$ Similarly we can do it for $2$-handle attachment. $\hfill \blacksquare$

\begin{theorem}\label{arf:2}
	Homotopy class of of $\Sigma(\eta \circ \Sigma \eta)$ in $\pi_5(\s^3)$ is non-trivial i.e. the Pontryagin manifold corresponding to $\Sigma(\eta \circ \Sigma \eta)$ in $\s^3$ is null cobordant.
\end{theorem}

\noindent \textbf{Proof}. First, we choose a regular value of \(\Sigma(\eta \circ \Sigma \eta)\), call it \( p \). Next, we choose the framing on \( p \) such that \( (\eta)^{-1}(p) \) has the twisted framing as discussed in the previous section. By the abundance of regular values, we may choose it in such a way that \( \Sigma(\eta \circ \Sigma \eta)^{-1}(p) = (\eta \circ \Sigma \eta)^{-1}(p) \), which is a circle bundle over the twisted circle \( \eta^{-1}(p) \). We can give it a framing coming from the framing on \( p \), making it an orientable \( 2 \)-manifold in \( \mathbb{S}^5 \). Hence, it must be a torus. Let \( \mathbb{T} \) be this torus. We will show that \( \arf(\mathbb{T}) = 1 \).  

\vspace{0.2cm}  

\noindent To see this, note that \( \mathbb{T} = (\eta \circ \Sigma \eta)^{-1}(p) \). Now, \( a = \Sigma \eta^{-1}(p) \) is framed, and \( \res(a) = 1 \) since it represents the generator of \( \pi_5(\mathbb{S}^4) \) via the correspondence. We can embed \( a \) in \( \mathbb{S}^5 \), and the residue class remains unchanged. Around \( a \), take a tubular neighborhood in \( \mathbb{S}^5 \), call it \( N \). Now, \( \nu: N \to a \) is a vector bundle, and \( \nu^{-1}(x) \) is homeomorphic to \( \mathbb{R}^4 \). In \( \mathbb{R}^4 \), we have a framed manifold corresponding to \( \eta \); call this \( b \). Then \( \res(b) = 1 \) since it corresponds to the generator of \( \pi_4(\mathbb{S}^3) \). The elements \( a \) and \( b \) can be visualized in the following picture:  

\[
	\includegraphics[width=0.45\textwidth]{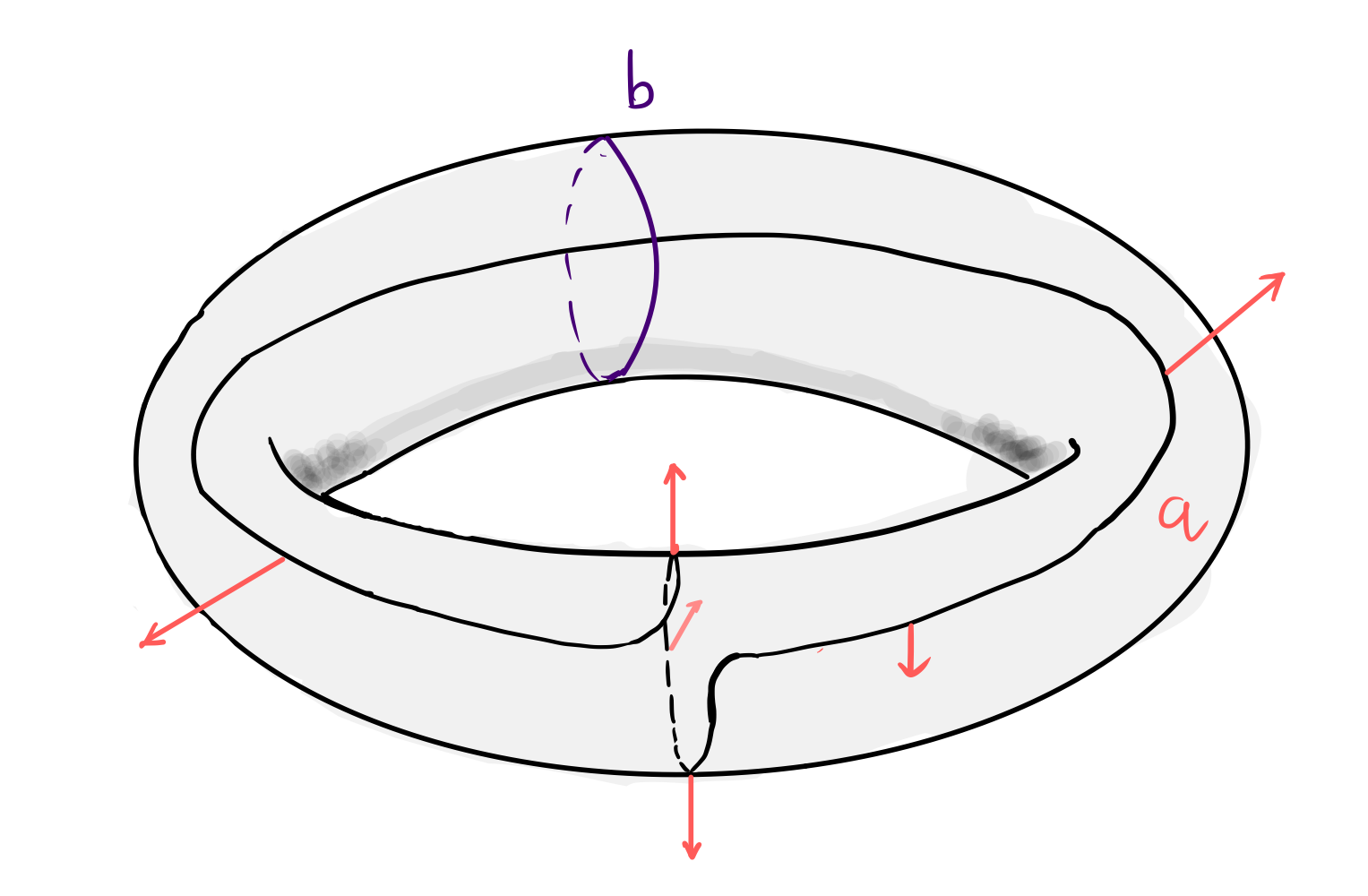}
\]  \begin{center}
     (This picture is taken from \cite{Hop1})
\end{center}

\noindent Note that \( a \) and \( b \) intersect at exactly one point and clearly represent generators of \( H_1(\mathbb{T}; \mathbb{Z}/2\mathbb{Z}) \). Thus,  
\[
\arf(\mathbb{T}) = q_{\mathbb{T}}(a) q_{\mathbb{T}}(b) = \res(a) \res(b) = 1,
\]
which completes the proof.

\section{Conclusion}

\noindent In this paper, we revisited Pontryagin’s original approach to computing the first two stable stems, emphasizing the role of framed submanifolds and their intersections. By carefully analyzing the framed bordism representatives and their corresponding residues, we recovered the classical computations of \(\pi_1^s\) and \(\pi_2^s\), highlighting the geometric intuition behind Pontryagin’s method. Our exposition clarifies how intersections encode secondary algebraic structures and how the framing conditions influence the calculations. This perspective not only reaffirms the validity of Pontryagin’s approach but also provides a foundation for extending similar techniques to higher-dimensional stable homotopy groups.  A similar analysis can be carried out for the third stable stem, as was done by Rokhlin. However, as can be observed, the difficulty of applying this method increases with the degree of the stable stem. Additionally, in this paper, we introduced certain residue classes and Arf invariants, which Pontryagin did not define using homology. His approach was to construct everything from first principles, while these modern constructions align with his work and provide a natural way to reinterpret his proof in a contemporary framework.

\end{document}